\documentclass[11pt]{article}
\usepackage{tipa}
\usepackage{latexsym,amssymb,amsfonts}
\usepackage{mathrsfs}
\usepackage{graphicx}
\usepackage{psfrag}
\usepackage{amsmath, amsbsy}
\usepackage{amsopn, amstext}
\usepackage{amsmath,multicol,amsopn}
\usepackage{latexsym,amssymb}
\usepackage{amsbsy, amstext,makeidx}

\def\sqr#1#2{{\vcenter{\vbox{\hrule height.#2pt
              \hbox{\vrule width.#2pt height#1pt \kern#1pt \vrule width.#2pt}
              \hrule height.#2pt}}}}
\def\signed #1{{\unskip\nobreak\hfil\penalty50
              \hskip2em\hbox{}\nobreak\hfil#1
              \parfillskip=0pt \finalhyphendemerits=0 \par}}
\def\endpf{\signed {$\sqr69$}}

\def\dbR{{\mathop{\rm l\negthinspace R}}}
\def\dbC{{\mathop{\rm l\negthinspace\negthinspace\negthinspace C}}}

\def\3n{\negthinspace \negthinspace \negthinspace }
\def\1n{\negthinspace }

\def\dbC{{\mathop{\rm l\negthinspace\negthinspace\negthinspace C}}}

\def\dbN{{\mathop{\rm l\negthinspace N}}}
\def\dbP{{\mathop{\rm l\negthinspace P}}}
\def\dbR{{\mathop{\rm l\negthinspace R}}}
\def\={\buildrel \triangle \over =}

\def\ds{\displaystyle}

\def\ns{\noalign{\ss}}

\def\d{\delta}
\def\e{\varepsilon}

\def\k{\kappa}

 \def\n{\nabla}

\def\f{\varphi}

\def\o{\omega}

\def\D{\Delta}

\def\Si{\Sigma}

\def\O{\Omega}

\def\cA{{\cal A}}

\def\cE{{\cal E}}

\def\cQ{{\cal Q}}
\def\cR{{\cal R}}

\def\cU{{\cal U}}

\def\cl{{\cal l}}
\def\no{\noindent}

\def\ss{\smallskip}
\def\ms{\medskip}
\def\bs{\bigskip}
\def\q{\quad}
\def\qq{\qquad}
\def\hb{\hbox}

\def\lan{\mathop{\langle}}
\def\ran{\mathop{\rangle}}

\def\esssup{\mathop{\rm esssup}}

\def\pa{\partial}

\def\wt{\widetilde}
\def\cd{\cdot}
\def\cds{\cdots}

\def\ae{\hbox{\rm a.e.{ }}}

\def\supp{\hbox{\rm supp$\,$}}
\def\span{\hbox{\rm span$\,$}}

\def\cl{\overline}
\def\co{\mathop{{\rm co}}}
\def\coh{\mathop{\overline{\rm co}}}

\def\|{\Big |}
\def\({\Big (}
\def\){\Big )}
\def\[{\Big[}
\def\]{\Big]}
\def\bde{\begin{definition}}
\def\ede{\end{definition}}
\def\be{\begin{equation}}
\def\bel{\begin{equation}\label}
\def\ee{\end{equation}}
\def\bt{\begin{theorem}}
\def\et{\end{theorem}}
\def\bc{\begin{corollary}}
\def\ec{\end{corollary}}
\def\bl{\begin{lemma}}
\def\el{\end{lemma}}
\def\bp{\begin{proposition}}
\def\ep{\end{proposition}}
\def\bas{\begin{assumption}}
\def\eas{\end{assumption}}
\def\br{\begin{remark}}
\def\er{\end{remark}}
\def\ba{\begin{array}}
\def\ea{\end{array}}
\def\ed{\end{document}}

\def\square#1{\vbox{\hrule\hbox{\vrule height#1%
     \kern#1\vrule}\hrule}}
\def\rectangle#1#2{\vbox{\hrule\hbox{\vrule height#1%
     \kern#2\vrule}\hrule}}

\font\tenbb=msbm10 \font\sevenbb=msbm7 \font\fivebb=msbm5

\newfam\bbfam
\scriptscriptfont\bbfam=\fivebb \textfont\bbfam=\tenbb
\scriptfont\bbfam=\sevenbb

\newtheorem{lemma}{Lemma}[section]
\newtheorem{remark}{Remark}[section]

\newtheorem{theorem}{Theorem}[section]
\newtheorem{corollary}{Corollary}[section]

\newtheorem{definition}{Definition}[section]
\newtheorem{proposition}{Proposition}[section]
\newtheorem{assumption}{Assumption}[section]

\makeatletter
   
   \@addtoreset{equation}{section}
\makeatother

\begin{document}

\title{\bf Finite Codimensional Controllability, and Optimal Control Problems with Endpoint State Constraints\thanks{This work is partially supported by the NSF of China
under grants   11871142 and  11471231,  by the
Fundamental Research Funds for the Central
Universities under grant 2412015BJ011, and by
PCSIRT under  grant IRT$\_$15R53.}}
\author{Xu Liu\thanks{Key Laboratory of Applied Statistics of MOE,  School of Mathematics and Statistics, Northeast Normal
University, Changchun 130024, China. E-mail
address: liux216@nenu.edu.cn.}\and Qi
L\"u\thanks{School of Mathematics,  Sichuan
University, Chengdu 610064,  China. E-mail
address: lu@scu.edu.cn.}\and Xu
Zhang\thanks{School of Mathematics, Sichuan
University, Chengdu 610064, China. E-mail
address: zhang\_xu@scu.edu.cn.}}

\date{}

\maketitle

\begin{abstract}
In this paper, motivated by the study of optimal
control problems for infinite dimensional
systems with endpoint state constraints, we
introduce the notion of finite  codimensional
(exact/approximate) controllability. Some
equivalent  criteria on the finite codimensional
controllability are presented. In particular,
the finite codimensional exact controllability
is reduced to deriving a G{\aa}rding type
inequality for the adjoint system, which is
new for many evolution equations. This
inequality can be verified for some concrete
problems (and hence applied to the corresponding
optimal control problems), say the wave
equations with both time and space dependent
potentials. Moreover, under some mild
assumptions, we show that the  finite
codimensional exact controllability of this sort
of wave equations is equivalent to the classical
geometric control condition.
\end{abstract}

\bs

\no{\bf Key Words}.  Finite  codimensional
controllability,
  finite codimensionality,   optimal control, endpoint state constraint, Pontryagin type
maximum principle.

\ms

\no{\bf AMS subject classifications}. 93B05, 49J20,
93B07, 49K20, 35Q93.

\section{Introduction}

It is well known that control theory was founded
by N. Wiener in 1948 (\cite{Wiener}). After
that, this theory was greatly extended to
various complicated setting and widely used in
sciences and technologies. Particularly, after
the  seminal works \cite{Bellman, Bellman1,
Kalman2,  Pontryagin, PC},  rapid development of
mathematical control theory (for both
deterministic and stochastic systems but this
paper will focus only on deterministic ones)
began in the 1960s (e.g., \cite{BLR, Clarke,  C,
Coron,  Crandall-Lions, Ekeland, Frankowska,
fur, Imanouilov, IK, 2, Lions, L, Robbiano, r,
Sussmann, xu, Zuazua, zu} and  rich references
cited therein). Usually, in terms of the
so-called state-space technique, people describe
the considered control system as a suitable
state equation.

Roughly speaking, ``control" means that one
hopes to change the dynamics of the involved
system, by means of a suitable way. In our
opinion, there are two (most, in some sense)
fundamental issues in control theory, i.e.,
feasibility and optimality, which we shall
explain more below.

The first fundamental issue is {\it
feasibility}, or in the terminology of control
theory, {\it controllability}, which means that,
one can find at least one way to achieve a goal.
More precisely, for simplicity, let us consider
the following controlled system governed by a
linear ordinary differential equation:
 \bel{ols1}
 \left\{
 \ba{ll}
 \ds y_t(t) =Ay(t)+Bu(t),\qq t>0,\\
  \ns
 y(0)=y_0.
 \ea\right.
 \ee
In (\ref{ols1}), $A\in \mathbb{R}^{n\times n}$,
$B\in \mathbb{R}^{n\times m}$ ($n,
m\in\mathbb{N}$), $y(\cd)$ is the state
variable, $u(\cd)$ is the control variable, and
$\mathbb{R}^n$ and $\mathbb{R}^m$ are the state
space and control space, respectively. The
system (\ref{ols1}) is called exactly
controllable (at  time $T>0$) if for any initial
state $y_0\in\mathbb{R}^n$ and any final state
$y_1\in\mathbb{R}^n$, there is a control
$u(\cd)\in L^2(0,T;\mathbb{R}^m)$ such that the
solution $y(\cd)$ to (\ref{ols1}) satisfies
\bel{eqx1} y(T)=y_1. \ee

The above definition of controllability can be
easily extended to abstract evolution equations.
In the general setting, it may happen that the
requirement (\ref{eqx1}) has to be relaxed in
one way or another. This leads to the
approximate controllability, null
controllability, partial controllability, and
finite codimensional  exact/approximate
controllability (to be introduced in this
paper), etc. Also, the above $B$ can be
unbounded for general controlled systems.

Clearly, the above controllability  problem can
be viewed as another equation problem, in which
both $y(\cd)$ and $u(\cd)$ are unknowns. Namely,
instead of viewing $u(\cd)$ as a control
variable, we may simply regard it as another
unknown variable. Nevertheless, the resulting
equation problem is definitely ill-posed.
Indeed, none of existence, uniqueness and
continuous dependence of this equation problem
is guaranteed. This is the main difficulty in
the study of many controllability problems (both
theoretically and numerically).

Controllability is strongly related to (or in
some situation, even equivalent to) other important issues
in control theory, say observability,
stabilization and so on.    One can find
numerous literatures on these topics (see
\cite{AMM, Barbu, BLR, BBZ, BZ, C, Coron,
Coron1, fur, Imanouilov, Ima-Yam, Kalman2, Le
Rousseau, LLTT, Lebeau, L, Robbiano, r,
Sussmann, xu, Zuazua, zu} and rich references
therein).

The second fundamental issue is {\it
optimality}, or in the terminology of control
theory, {\it optimal control}, which means that
people are expected to find the best way, in
some sense, to achieve their goal. As an
example, we fix $y_0,y_1\in\dbR^n$. It is easy
to see that, if there exists a control $u(\cd)$
such that the solution $y(\cd)$ to (\ref{ols1})
satisfies (\ref{eqx1}), then very often one may
find another control verifying the same
conditions. Naturally, one hopes to find the
``best" control fulfilling these conditions. To
be more precisely, we fix a suitable function
$f^0(\cd,\cd,\cd):[0,T]\times\mathbb{R}^n\times\mathbb{R}^m\to\mathbb{R}$,
and denote by $\mathcal{U}_{ad}$ the set of
controls $u(\cd)\in L^2(0,T;\mathbb{R}^m)$ so
that the solution $y(\cd)$ to (\ref{ols1})
satisfies (\ref{eqx1}) and
$f^0(\cd,y(\cd),u(\cd))\in L^1(0,T)$. A typical
optimal control problem for the system
(\ref{ols1}) is to find a $\bar
u(\cd)\in\mathcal{U}_{ad}$, called an optimal
control, which minimizes the following
functional over $\mathcal{U}_{ad}$:
 $$
 J(u(\cd))=\int_0^Tf^0(t,y(t),u(t))dt.
 $$

The above formulation of optimal control problem
can  be easily extended to more general setting,
in particular for the case that the controls
take values in a more general set (instead of
$\mathbb{R}^m$), say any metric space, which
does not need to enjoy any linearity or
convexity structure.

Optimal control problems are strongly related to
the classical calculus of variations and
optimization theory. Nevertheless, since the
control set may be quite general, the classical
variation technique cannot be applied to optimal
control problems directly, especially in the
case that the state space is infinite
dimensional. Various optimal control problems
are extensively studied in the literatures
(e.g., \cite{Bellman, Bellman1, Clarke,
Crandall-Lions, Ekeland, Frankowska, IK, LY,
Lions, Pontryagin, PC, WWXZ}  and rich
references cited therein).

Clearly, the study of controllability problems
is a basis to investigate further optimal
control problems. Indeed, the usual nonempty
assumption on the set of feasible/admissible
control-state pairs (for optimal control
problems) is actually a controllability
condition. Nevertheless, in the previous
literatures, it seems that the studies of
controllability and optimal control problems are
almost independent. Two typical exceptions that
we know are the following:
\begin{itemize}

   \item[1)] In \cite{Ima-Yam}, some techniques from optimal control theory are employed to derive the observability estimate and null controllability for parabolic type equations.

   \item[2)] In \cite{WWXZ}, some techniques developed in the study of controllability and observability problems are adopted to solve several time optimal control problems.
\end{itemize}

In our opinion, now it is the time to solve
controllability and optimal control problems as
a whole, at least in some sense and to some
extend, though they are two different control
issues. This is by no means an easy task.
Actually,  for many concrete problems, it is
highly nontrivial to verify the above mentioned
assumption that the set of feasible/admissible
control-state pairs is nonempty.

The main purpose of this paper is to provide a
new link  between controllability and optimal
control problems in infinite dimensions. Our
work is motivated by the study of optimal
control problems for abstract evolution
equations with endpoint state constraints. In
\cite{Fattorini, ly, LY}, in order to guarantee
the nontriviality of Lagrange type multipliers in the
corresponding Pontryagin type maximum principle,
a finite codimensionality condition is
introduced. However, it is usually very
difficult to verify this condition directly
except for some  very special cases. Because of
this,  we shall reformulate this condition as a
class of new controllability notion, i.e. finite
codimensional exact controllability.

A key contribution in this work is to reduce
further the above mentioned finite codimensional
exact controllability to a suitable {\it a
priori} estimate for the underlying adjoint
system (see the estimate  (\ref{5.21-eq1}) for
the equation (\ref{215})). We remark that, in
some sense, the inequality (\ref{5.21-eq1}) can
be regarded as a G{\aa}rding type inequality,
which concerns the lower bound of  a bilinear
form induced by a linear elliptic
(pseudo-)differential operator. To see this, let
us recall below the classical G{\aa}rding
inequality (e.g., \cite[Section 5 of Chapter
X]{Hormander} for more details and more general
results). Let   $\Omega$ be a bounded domain in
$\dbR^n$ with   a  smooth boundary
$\partial\Omega$, and let $L$ be a uniformly
linear elliptic differential operator of order
$2k$ (for some positive integer $k$) with smooth
coefficients, i.e., there exists a constant
$s_0> 0$ such that
$$
\ell(x,\xi)\geq s_0|\xi |^{2k},\qq \forall\;
(x,\xi )\in \Omega\times\dbR^{n},
$$
where $\ell$ is the symbol of $L$. Then there
exist two constants $C_0>0$ and $C_1\geq 0$ such
that
\begin{equation}\label{9.7-eq1}
C_0|v|_{H^{k}(\Omega )}^{2}\leq \langle L
v,v\rangle_{H^{-k}(\Omega),H_{0}^{k}(\Omega)}
+C_1|v|_{L^{2}(\Omega)}^{2} ,\qq \forall\;v\in
H_{0}^{k}(\Omega).
\end{equation}
Clearly, both (\ref{5.21-eq1})  and
\eqref{9.7-eq1}  have an extra term, that is,
$|G\phi_T|_X$ and $|v|_{L^{2}(\Omega)}^{2}$,
respectively. It is easy to observe that these
two terms are accordingly compact with respect to the ones in
the left hand sides of  the corresponding
estimates. Hence, we may call (\ref{5.21-eq1}) a
G{\aa}rding type inequality for the evolution
equation (\ref{215}). This inequality can be
verified for some concrete problems, say the
wave equations with both time and space
dependent potentials (see Subsection
\ref{sub6.1}). Though the later result (which
seems not available in the previous literatures)
might be known for some experts in the field of
micro-local analysis, nobody knows how to use
it. Interestingly, in this work we shall give
its application in optimal control problems.
Moreover, under some mild assumptions, we shall
show that the  finite codimensional exact
controllability of this sort of wave equations
is equivalent to the celebrated geometric
control condition  (introduced in the papers
\cite{BLR1, BLR}) for the classical wave equation.

In this work, in order to present the key idea
in  the simplest way, we shall not pursue the
full technical generality. It deserves
mentioning that the method and technique
developed in this paper can be employed to
handle many other problems. Especially, our
finite codimensionality technique can be applied
to solve some interesting problems in
optimization, calculus of variations and
stochastic control, and even gives new results
for some finite dimensional optimal control
problems under state constraints (see our
forthcoming paper \cite{LLZ} for more details).

The rest of this paper is organized as follows.
Section \ref{sss3} is of preliminary nature, in
which we present some notations, notions and
simple results. In Section \ref{sec-ex}, some
equivalent criteria for the finite codimensional
exact controllability are given. Section
\ref{sec-ap} is devoted to a characterization of
the finite codimensional approximate
controllability. In Section \ref{sec-op},  the
finite codimensional controllability is applied
to study some optimal control problems  with
state  constraints. In Section \ref{sec-exam},
two examples    are given.  Finally, in
Appendix, we prove a technical result used in
this paper.

\section{Notations, notions and some preliminary results}\label{sss3}

  To begin with, we  introduce some notations.
Let   $Y$  and $U$ be two  reflexive Banach
spaces. For a Banach space $Z$,  denote by
$\mathcal{L}(Z;Y)$ the set of all bounded linear
operators from $Z$ to $Y$,  and  write it
$\mathcal{L}(Y)$  for short when $Z=Y$. For any
operator $P\in \mathcal{L}(Z;Y)$, write $P^*$
for its adjoint operator. Denote by $Y'$ the
dual space of $Y$; by $\overline{D}$ the closure
of a subset $D$ of $Y$;   by span$D$ the closed
subspace spanned by $D$;  and by $\coh D$ the
convex closed hull  of $D$.   For two subsets
$D_1$ and $D_2$ of $Y$, set $D_1-D_2=\big\{y\in
Y\ \big|\  y=y_1-y_2\mbox{  for some  }y_1\in
D_1 \mbox{ and }y_2\in D_2\big\}$ and
$D_1\setminus D_2=\big\{y\in Y\ \big|\  y\in
D_1\mbox{  and  }y\notin D_2\big\}$.
 Let $T>0$, $p\in (1, \infty]$,
$\mathcal{U}_p=L^p(0, T; U)$ and  $i$  be the
usual imaginary unit.

\smallskip

Consider the following  linear  control  system:
\begin{eqnarray}\label{11}
\left\{
\begin{array}{ll}
\ds y_t(t)=Ay(t)+F(t)y(t)+B(t)u(t),  &\quad t\in (0, T],\\
\ns\ds y(0)=y_0,
\end{array}\right.
\end{eqnarray}
where  $u(\cdot)\in\mathcal{U}_p$ is the control
variable and $y(\cdot)$ is the state variable,
$y_0\in Y$,  $A: \mathcal{D}(A)\subseteq
Y\rightarrow Y$  generates a $C_0$-semigroup on
$Y$, $F(\cdot)\in L^\infty(0, T;
\mathcal{L}(Y))$, and $B(\cdot)\in L^\infty(0,
T; \mathcal{L}(U; Y))$. Obviously,  $A+F(\cdot)$
generates an evolution operator
$\mathcal{S}(\cdot, \cdot)$ on $Y$. For any
$T>0$,   $y_0\in Y$ and $u(\cdot)\in
\mathcal{U}_p$, (\ref{11}) admits a mild
solution $y(\cdot)=y(\cdot; y_0, u(\cdot))\in
C([0, T]; Y)$, and $$ y(t)=y(t; y_0,
u(\cdot))=\mathcal{S}(t,
0)y_0+\displaystyle\int^t_0 \mathcal{S}(t,
s)B(s)u(s)ds, \q \forall\ t\in [0, T]. $$
Define the reachable   set  $\mathcal{R}(T;
y_0)$ of (\ref{11}) at time $T$ with the initial
value $y_0$ as follows:
$$
\mathcal{R}(T; y_0)=\Big\{y(T; y_0, u(\cdot))\in
Y\ \Big|\ y(\cdot) \mbox{ is the mild solution
to (\ref{11}) with some } u(\cdot)\in
\mathcal{U}_p\Big\}.
$$

Next, let us recall the notions of  finite
codimensional subspace  and finite
codimensionality (e.g., \cite{2}).

\begin{definition}
A  linear subspace $Y_0$ of  $Y$  is called
finite codimensional, if there exist an
$m\in\dbN$ and linearly independent $y_1, y_2,
\cdots, y_m\in Y\setminus Y_0$ such that $
\mbox{span}\big\{ Y_0, y_1, y_2,  \cdots, y_m
\big\}=Y. $
\end{definition}
\begin{definition}\label{10.25-def1}
A subset $D$ of  $Y$ is called finite
codimensional  in $Y$,   if

\smallskip

\noindent ${\bf (H_1)}$ There exists a
 $y_0\in\coh D$,   such that
$\span\big\{D-y_0\big\}$  is a finite
codimensional subspace of $Y$; and

\smallskip

\noindent  $\bf  (H_2)$ $\coh \big(D-y_0\big)$
has at least  an interior point in this
subspace.
\end{definition}

Now, we introduce two notions of finite
codimensional controllability.
\begin{definition}\label{d1}
The system $(\ref{11})$ is called finite
codimensional  exactly $($resp.,
approximately$)$ controllable at time $T$,  if
$\mathcal{R}(T;  0))$ $($resp.,
$\overline{\mathcal{R}(T; 0)}$ is  a finite
codimensional subspace of $Y$.
\end{definition}
\begin{remark}
Recall that  $(\ref{11})$ is exactly $($resp.,
approximately$)$ controllable at time $T$, if
$\mathcal{R}(T;  0)=Y$ $($resp.,
$\overline{\mathcal{R}(T; 0)}=Y)$. Therefore,
the finite codimensional exact $($resp.,
approximate$)$ controllability defined in
Definition $\ref{d1}$ is clearly weaker than the
usual exact $($resp., approximate$)$
controllability for linear systems.
\end{remark}
\begin{remark}
In general,   the  finite codimensional exact
controllability   cannot be reduced to the usual
exact controllability problem. Indeed, this is
possible only for the special case that $A+F(t)$
in $(\ref{11})$ has an invariant subspace,
which is  finite codimensional in $Y$ and
independent of $t\in  [0, T]$.
\end{remark}

As mentioned before, the notion of  finite
codimensional controllability is motivated by
the study  of  some optimal control problems for
infinite dimensional systems with endpoint state
constraints. It is well known that Pontryagin's
maximum principle is one of the milestones in
optimal control theory.   As a necessary
condition of optimal controls, for very general
finite dimensional systems, Pontryagin type
maximum principle was established in
\cite{PC}. Nevertheless,   surprisingly, it
fails for infinite dimensional systems if there
is  no further assumption (see \cite{e}). This
leads to that for quite a long time, Pontryagin
type maximum principle had been studied only for
evolution equations without terminal state
constraints. Until 1980s, by assuming the
finite  codimensionality of  some subset in
state spaces, Pontryagin type maximum principles
on optimal control problems for infinite
dimensional systems with endpoint  constraints
and  general control domains  were established
in  \cite{Fattorini, ly, LY}. In the following,
we present an optimal control problem with
state constraints and recall how  to use the
finite codimentionality   in deriving Pontryagin
type  maximum principle.

\smallskip

Consider the following
 evolution equation on
$Y$:
\begin{equation}\label{61}
\ds y_t(t)=Ay(t)+f(t, y(t), u(t)),  \q  t\in (0,
T],
\end{equation}
where $u(\cdot)$ is the control variable and
$y(\cdot)$ is the state variable.   Assume that
$f: [0, T]\times Y\times U\rightarrow Y$
satisfies certain conditions (to be given
later), such that for any $y(0)\in Y$ and
$u(\cdot)\in\mathcal{U}_p$, (\ref{61}) admits a
mild solution $y(\cdot)=y(\cdot;  y(0),
u(\cdot))\in C([0, T];  Y)$. Also, let
$\widetilde U$  be a nonempty subset of  $U$,
and
 $S$ be a closed and convex subset of $Y\times Y$.  Put
\begin{eqnarray*}
\begin{array}{ll}
\ds\mathcal{U}[0, T]=\Big\{
u(\cdot)\in\mathcal{U}_p\ \Big|\ u: (0,
T)\rightarrow \widetilde{U}\mbox{  is
measurable}  \Big\},& \\\ns\ds
\ds\mathcal{U}_{ad}=\Big\{
u(\cdot)\in\mathcal{U}[0, T]\ \Big|\ \mbox{the
mild  solution }y(\cdot)\mbox{ to }(\ref{61})
\mbox{  satisfies  }\big (y(0),y(T)\big)\in S
\Big\},&\\\ns\ds \mathcal{A}_{ad}=\Big\{
\big(u(\cdot),
y(\cdot)\big)\in\mathcal{U}_{ad}\times C([0,T];
Y)\ \Big|\  y\mbox{ is the  mild solution to
}(\ref{61})\Big\}.&
\end{array}
\end{eqnarray*}
Write
$$
J\big(u(\cdot),
y(\cdot)\big)=\displaystyle\int^T_0 f^0(t, y(t),
u(t))dt,
$$
where  $f^0:  [0, T]\times Y\times U\rightarrow
\dbR$ satisfies certain conditions (to be
specified in the sequel),  such that for any
$u(\cdot)\in \mathcal{U}[0, T]$, $y(0)\in Y$ and
the corresponding mild solution $y(\cdot)$ to
(\ref{61}),     $f^0\big(\cdot, y(\cdot),
u(\cdot)\big)\in L^1(0, T)$.

Consider the following optimal control problem
for  the system (\ref{61}):

\smallskip

\noindent {\bf (P)}\quad Find  a  pair
$(\overline{u}(\cdot), \overline{y}(\cdot))
 \in \mathcal{A}_{ad}$, such that
$J(\overline{u}(\cdot),
\overline{y}(\cdot))=\inf\limits_{(u(\cdot),
y(\cdot))\in \mathcal{A}_{ad} } J(u(\cdot),
y(\cdot)). $

\noindent Such a $(\overline{u}(\cdot),
\overline{y}(\cdot))$  is called an optimal
pair. As a  necessary condition for optimal
pairs,   Pontryagin type maximum principle  is
stated as follows.

\smallskip

\noindent{\bf Pontryagin type maximum
principle}: Assume that $(\bar u(\cd),\bar
y(\cd))$ is an optimal pair.   Then there exists
a pair $(\psi^0,\psi(\cd))\in\dbR\times C([0,T];
Y')$,
 such that
\begin{equation}\label{0llz}
\big(\psi^0,  \psi(\cdot)\big)\neq \big(0,
0\big),
\end{equation}
\begin{equation}\label{6.5-eq1}
\psi_t(t)= -A^*\psi(t) - f_y(t,\bar y(t), \bar
u(t))^*\psi(t)-\psi^0f^0_y(t,\bar y(t),\bar
u(t)),\q\mbox{a.e. } t\in (0,T),
\end{equation}
\begin{equation}\label{6.5-eq2}
\langle\psi(0), y^0-\bar y(0)\rangle_{Y',
Y}-\langle\psi(T), y^1-\bar y(T)\rangle_{Y',
Y}\le0,\q\forall\   (y^0, y^1)\in S,
\end{equation}
\begin{equation}\label{6.5-eq3}
H\big(t,\bar y(t),\bar u(t),\psi^0,\psi(t)\big)=
\max_{u\in \widetilde{U}}H\big(t,\bar y(t),u,
\psi^0,\psi(t)\big), \q \mbox{a.e. } t\in(0,T),
\end{equation}
where  $A^*$ is the  adjoint operator of  $A$,
and
$$
H(t,y,u,\psi^0,\psi)\ds=\psi^0f^0(t,y,u)
+\langle\psi,f(t,y,u)\rangle_{Y', Y},
$$
$$\q\q\q\q\q\q\q\q\q\q\q \ds\forall\    (t,y,u,\psi^0,\psi)\in[0,T]\times
Y\times \widetilde{U}\times\dbR\times Y'.
$$

\smallskip

\eqref{0llz} is key  in Pontryagin type maximum
principle. Indeed, if it fails, then $\psi^0=0$
and  $\psi(t)=0$ for all $t\in [0,T]$. Hence,
\eqref{6.5-eq2} and \eqref{6.5-eq3} are trivial,
since they are  then simply  ``$0\leq 0$" and
``$0=0$", respectively.

\smallskip

In order to ensure \eqref{0llz},  the finite
codimensionality of a suitable  set was
introduced.  More precisely, consider the
following system:
\begin{eqnarray}\label{10.25-eq4}
\left\{
\begin{array}{ll}
\ds\xi_t(t)= A\xi(t) + f_y(t,\bar y(t), \bar
u(t))\xi(t)+f(t,\bar y(t),u(t))-f(t,\bar
y(t),\bar u(t)),  &\quad   t\in (0,T],\\
\ns\ds\xi(0)=0,
\end{array}
\right.
\end{eqnarray}
and the homogenous equation:
\begin{eqnarray}\label{10.25-eq5}
\left\{
\begin{array}{ll}
\ds \eta_t(t)= A\eta(t) + f_y(t,\bar
y(t), \bar u(t))\eta(t), &\quad   t\in (0,T],\\
\ns\ds\eta(0)=y^0,
\end{array}\right.
\end{eqnarray}
for $y^0\in Y$.  Put
$$
\cR=\Big\{ \xi(T) \in Y\  \Big|\  \xi(\cd)
\mbox{ is the mild solution to }
\eqref{10.25-eq4}\mbox{  with some }
u(\cd)\in\cU[0,T]  \Big\}
$$
and
$$
\cQ=\Big\{ y^1-\eta(T)\in Y\ \Big|\   \eta(\cd)
\mbox{ is the mild solution to }
\eqref{10.25-eq5}\mbox{ and }  (y^0,  y^1)\in S
\Big\},
$$
and introduce the condition:
$$
{\bf (H)} \q  \cR\!-\!\cQ \mbox{ is finite
codimensional  in } Y.
$$

\noindent  It was  proved in \cite[Chapter
4]{2} that,
% the following necessary condition on  the optimal  pair $(\bar u(\cd),\bar y(\cd))$
%was proved in \cite{2}.
%\begin{proposition}\label{10.25-th1}
if the condition ${\bf (H)}$ holds,   the
optimal pair $(\bar u(\cd),\bar y(\cd))$ in the
optimal  control  problem {\bf (P)} satisfies
Pontryagin type maximum  principle, i.e.,
$(\ref{0llz})$-$(\ref{6.5-eq3})$  hold.
%\end{proposition}
%

\smallskip

The proof    is  based on  the following known
result.
\begin{lemma}\label{LM}
$(\cite[\mbox{Lemma } 3.6  \mbox{ on Page  }
142]{2})$ If  $M$ is finite codimensional in
$Y$,  then  for any\linebreak
$\{f_j\}_{j=1}^\infty\subseteq Y'$ satisfying
the following two conditions:

\medskip

\noindent $(1)$ $|f_j|_{Y'}\geq \delta$ for a
positive constant  $\delta$  and $f_j\rightarrow
f$ weakly$^*$ in $Y'$, as $j\rightarrow \infty$;
and

\medskip

\noindent $(2)$ There exist positive constants
$\epsilon_j$, such that
$\lim\limits_{j\rightarrow \infty} \epsilon_j=0$
and $\lan f_j, x\ran_{Y',Y}\geq -\epsilon_j$,
 $\forall\ x\in M$,

\noindent it holds that $f\neq 0$.
\end{lemma}

Lemma \ref{LM}  means that, under some mild
assumptions, the finite codimensionality  on $M$
is sufficient to guarantee the weak limit point
of a sequence  to be nonzero in an infinite
dimensional  space. This is the reason why
(\ref{0llz}) holds in Pontryagin type maximum
principle. On the other hand, this condition is
also necessary, at least when $Y$ is a Hilbert
space and $M$ is a linear closed subspace. In
fact, we have the following result.
\begin{proposition}\label{PP1}
Suppose that  $M$ is a linear closed subspace of
a Hilbert space $Y$. Then $M$ is  finite
codimensional, if and only if for any
$\{f_j\}_{j=1}^\infty\subseteq Y$ satisfying the
conditions $(1)$-$(2)$ in Lemma $\ref{LM}$, it
holds that $f\neq 0$.
\end{proposition}

\noindent{\bf Proof. } By Lemma \ref{LM},  we
only need  to prove  the sufficiency. If $M$ is
not finite codimensional,  then there exists a
subspace $Y_0=\mbox{span}\{e_1, e_2, \cdots\}$
of $Y$, such that $M\oplus Y_0=Y$. Also,
$|e_j|_Y=1$ for any $j\in\dbN $ and
$\{e_j\}_{j=1}^\infty$ is pairwise orthogonal.
Choose $f_j=e_j$, $\delta=1$  and
$\epsilon_j=1/j$. Then $|f_j|_Y=1$,  $(f_j,
x)_Y=0\geq -1/j$ for any $x\in M$  and
$\lim\limits_{j\rightarrow \infty}
\epsilon_j=0$. Notice that
$$e_j\rightarrow \hat{e} \mbox{ weakly in }Y\
\mbox{ with }\ |\hat{e}|_Y\leq 1.$$ Since
$\sum\limits_{j=1}^{\infty} (e_j, x)_Y<\infty$
for  any $x\in Y$,
$\lim\limits_{j\rightarrow\infty}(e_j, x)_Y=0$,
which  implies that $\hat{e}=0$. This
contradicts the assumptions on
$\{f_j\}_{j=1}^\infty$ and therefore,  $M$ is
finite codimensional in $Y$.
\endpf

\begin{remark}
 Proposition $ \ref{PP1}$   indicates that the finite codimensionality   seems closely
related to the  weak  convergence method and
existence of nontrivial  solutions for partial
differential equations.  We shall study  its
applications in this respect  in a future work.
\end{remark}

Usually, unless $\cQ$ is finite codimensional in
$Y$,  it is quite difficult to verify  the
condition {\bf (H)} directly,   even for some
simple linear systems. For example, if
$S=\{y^0\}\times B_1$ for a given   $y^0\in Y$
and the unit ball $B_1$ of $Y$, then {\bf (H)}
holds trivially. But when $S=\{(y^0, y^1)\}$
with $y^0, y^1\in Y$, it seems not easy to check
this condition,  since the set $\cR$ is the
reachable set of some system with control
constraints. As we mentioned before, the
motivation of this paper is to introduce a new
method to verify the finite codimensionality
condition appeared in optimal control  problems.
A little more precisely, first, the condition
{\bf (H)} is reduced to a finite codimensional
exact controllability problem,  as introduced in
Definition \ref{d1}. Then,  by a duality
argument,  such a controllability problem is
transformed into a suitable {\it a priori}
estimate, called weak observability estimate
(compared to the usual observability estimate),
for its adjoint system, which is more easily
verified or proved false, at least for some
nontrivial examples  (see Propositions
\ref{llz11} and \ref{pllz}).

\section{Finite codimensional  exact
controllability}\label{sec-ex}

In this section,   some equivalent results on
finite codimensional exact  controllability are
established. First,  consider the following
linear control system:
\begin{eqnarray}\label{llz51}
\left\{
\begin{array}{ll}
\ds y_t(t)=Ay(t)+F(t)y(t)+B(t)u(t),  &\quad t\in (0, T],\\
\ns\ds y(0)=0,
\end{array}\right.
\end{eqnarray}
where  $A$, $F(\cdot)$ and $B(\cdot)$ are the
same as those in (\ref{11}). Assume that
$$
{\bf (A)}\q  \widetilde{\mathcal{U}} \mbox{ is
a nonempty bounded subset of  }
 \mathcal{U}_p \mbox{ and  }   \overline{\co}\ \wt \cU  \mbox{  has   at least
an  interior  point}.
$$
Set
\begin{equation}\label{llz52}
M=\Big\{ y(T)\in Y\ \Big|\ y\mbox{ is the mild
solution to } (\ref{llz51})\mbox{ with  some  }
u(\cdot)\in\widetilde{\mathcal{U}}  \Big\}.
\end{equation}
Then it is easy to check that
\begin{equation}\label{llz53}\left\{
\begin{array}{ll}
\ds\mbox{span}M=\overline{\Big\{ y(T)\in Y\ \Big|\ y\mbox{ is the mild solution to } (\ref{llz51})\mbox{ with some } u(\cdot)\in\mathcal{U}_p  \Big\}},&\\
\ns\ds \overline{\mbox{co}} M=\Big\{ y(T)\in Y\
\Big|\ y\mbox{ is the   mild solution to }
(\ref{llz51})\mbox{  with some  }
u(\cdot)\in\overline{\co}\ \wt \cU  \Big\}.&
\end{array}
\right.
\end{equation}
Also, we recall a known result on finite
codimensional subspace.

\begin{lemma}\label{lm1}
 $(\cite[\mbox{Proposition  }3.2\mbox{ on  Page }138]{2})$
Assume that $Y_0$ is a linear  subspace of $Y$.
Then  $Y_0$ is finite codimensional in $Y$, if
and only if there exist finitely many bounded
linear functionals $\{f_j\}_{j=1}^m\subseteq
Y'$,  such that
$Y_0=\bigcap\limits_{j=1}^m\ker\{f_j\}$.
\end{lemma}

The first result of this section is stated as
follows.
\begin{theorem}\label{t2}
Suppose that ${\bf (A)}$ holds.  Then the
following two assertions are equivalent:

\medskip

\noindent  {\bf (1) }  The system
$(\ref{llz51})$ is finite  codimensional exactly
controllable in $Y$.

\medskip

\noindent  {\bf (2) } The set $M$ $($in
$(\ref{llz52}))$ is finite codimensional in $Y$.

\end{theorem}

\noindent {\bf Proof.}  Without loss of
generality, we assume that $0$ is an interior
point of $\overline{\co}\  \wt \cU$. Otherwise,
if $u_0\neq 0$  is an interior point of
$\overline{\co}\  \wt \cU$, it suffices to
replace  $\wt\cU$ and $M$,  respectively, by
$\wt \cU-u_0$ and $M-y(T;  u_0)$ with $y(\cdot;
u_0)$ being   the mild solution to (\ref{llz51})
associated to $u=u_0$. Hence, there is an
$r_0>0$, such that $ \big\{ u(\cdot)\in \cU_p\
\big|\ |u|_{\cU_p}\leq r_0
\big\}\subseteq\overline{\co}\ \wt \cU\
\mbox{and}\  0\in M. $

\smallskip

First, we prove that  ${\bf (1)}$ implies ${\bf
(2)}$. For any  $n\in\dbN$, set
$$
N_n=\Big\{ y(T)\in Y\ \Big|\ y \mbox{ is the
mild solution to }(\ref{llz51})\mbox{   with
some }      u(\cdot)\in\mathcal{U}_p \mbox{
satisfying } |u|_{\mathcal{U}_p}\leq nr_0
\Big\}.
$$
Then $N_1\subseteq \overline{\co}  M$ and
$\bigcup\limits_{n\in\dbN} N_n =\mathcal{R}(T;
0)$. By  ${\bf (1)}$ and  (\ref{llz53}),
$\mathcal{R}(T; 0)=\overline{\mathcal{R}(T;
0)}=\mbox{span}  M$ is a finite   codimensional
subspace of $Y$. Also, by the Baire category
theorem,  there exists an $\tilde n\in\dbN$,
such that $N_{\tilde n}=\overline{N_{\tilde n}}$
has at least an interior point $\tilde  y \in
{\rm span}M$. Then $\ds\frac{\tilde y}{\tilde
n}$ is an interior point of $\overline{\co} M$
in span$M$. Hence, by Definition
\ref{10.25-def1}, $M$ is  finite codimensional
in $Y$.

\smallskip

On the other hand, we prove  that ${\bf (2)}$
implies ${\bf (1)}$.   Notice that
$\overline{\co}M\subseteq \mathcal{R}(T; 0)$.
By
 $({\bf H_2})$ in Definition
\ref{10.25-def1}, $\overline{\mbox{co}} M$ has
at  least  an  interior point in  the subspace
span$M=\overline{\mathcal{R}(T; 0)}$. Hence,
$\mathcal{R}(T; 0)$ also has an  interior  point
in $\overline{\mathcal{R}(T; 0)}$. Since
$\mathcal{R}(T; 0)$ and
$\overline{\mathcal{R}(T; 0)}$ are two linear
subspaces of $Y$  and $\mathcal{R}(T; 0)$ is
dense in $\overline{\mathcal{R}(T; 0)}$, it
follows that $\mathcal{R}(T;
0)=\overline{\mathcal{R}(T; 0)}$.  Also,  by
$({\bf H_1})$ in Definition \ref{10.25-def1},
$\mathcal{R}(T; 0)$=span$M$  is finite
codimensional  in  $Y$.  Hence, ${\bf (1)}$
holds.
\endpf

\medskip

Next, by a duality technique,  we prove that the
finite  codimensional exact controllability of
(\ref{llz51}) is  equivalent to  a suitable
observability estimate for the following
equation (or adjoint system):
\begin{eqnarray}\label{215}
\left\{
\begin{array}{ll}
\ds\phi_t(t)=-A^*\phi(t)-F(t)^*\phi(t),
&\quad t\in (0, T],\\
\ns\ds\phi(T)=\phi_T,
\end{array}\right.
\end{eqnarray}
where $\phi_T\in Y'$.  Set  $p'=p/(p-1)$ for
$p\in (1, \infty)$, and $p'=1$ for $p=\infty$.
In what follows, $C$ is used to denote a generic
positive constant,  which may change from line
to line in the sequel.

\medskip

The second  result  of this section is as
follows.
\begin{theorem}\label{t2.1}
The following  two assertions are equivalent:

\medskip

\noindent  {\bf (1) }  The system
$(\ref{llz51})$ is finite  codimensional exactly
controllable in $Y$.

\medskip

\noindent  {\bf (2) } There  exists a finite
codimensional subspace $\wt Y$ of  $Y'$,  such
that  any  solution $\phi$ to  $(\ref{215})$
satisfies
\begin{equation}\label{221}
|\phi_T|_{Y'}\leq C |B(\cdot)^*\phi|_{L^{p'}(0,
T; U')},\quad\quad\forall\,  \phi_T\in \wt Y.
\end{equation}

\end{theorem}

\noindent  {\bf Proof.  }  First,  we  prove
that ${\bf (2)}$ implies ${\bf (1)}$. The  proof
is divided into four parts.

\smallskip

\noindent{\bf Step 1.} In this step, we prove
that the following subspace $Y_1$ of $Y$  is
finite dimensional:
$$
Y_1=\big\{ x\in Y\ \big| \ \langle  f, x
\rangle_{Y', Y}= 0,\ \forall f\in \wt Y \big\},
$$
where  $\wt Y$ is the subspace given  in  {\bf
(2)}.

Let the codimension of $\wt Y$ be $k_1$. If
$Y_1$ is an infinite dimensional space, then
there is a linear subspace $Y_1^0$ of $Y_1$,
whose
 dimension is $k_1+1$.   Let
$\{x_1,\cds,x_{k_1+1}\}\subseteq Y$ be a basis
of $Y_1^0$. By the  Hahn-Banach theorem, one can
find $\{f_1,\cds,   f_{k_1+1}\}\subseteq Y'$,
such that for $1\leq   k,   j\leq k_1+1$,
\begin{eqnarray*}
\langle   f_k, x_j\rangle_{Y', Y}= \left\{
\begin{array}{ll}\ds
1 &\mbox{ for } k=j,\\
\ns\ds 0 &\mbox{ for } k\neq j.
\end{array}
\right.
\end{eqnarray*}
It follows that $\{f_1,\cds, f_{k_1+1}\}$ are
linearly independent in  $Y'$. Hence,  the
dimension of the subspace ${\rm
span}\big\{f_1,\cds,  f_{k_1+1}\big\}$ is
$k_1+1$. Since for any $j=1, \cdots, k_1+1$,
$x_j\in  Y_1$ and $\langle  f_j, x_j\rangle_{Y',
Y}=1\neq 0$,  by the definition of $Y_1$,  we
get that $ f_j\notin \wt Y$ $(j=1, \cdots,
k_1+1)$. This contradicts the fact that the
codimension of $\wt Y$ is $k_1$.  Hence,  $Y_1$
is finite  dimensional and
  denote by
$k_2$ its dimension.

\smallskip

\noindent {\bf Step 2.} In this step, we prove
that for any $y_T\in  Y$,   there is a control
$u\in \mathcal{U}_p$,   such that the
corresponding solution $y(\cdot;  u)$  to
(\ref{llz51}) satisfies
\begin{equation}\label{00llz00}
y(T; u)-y_T\in Y_1.
\end{equation}

Write $ K=\Big\{B(\cdot)^*\phi\in L^{p'}(0, T;
U')\ \Big|\ \phi \mbox{ is the mild solution to
(\ref{215}) with $\phi_T\in \wt Y$} \Big\} $ and
define a linear functional $\ell$ on $K$ as  $
\ell\big(B(\cdot)^*\phi\big)=\langle \phi_T, y_T
\rangle_{Y',   Y}. $
It follows from (\ref{221}) that $\ell$ is a
bounded linear functional on $K$. Moreover,
$|\ell|_{\mathcal{L}(K; \dbR)}\leq C|y_T|_Y$.
Hence, by the Hahn-Banach theorem, $\ell$ can be
extended   to be a bounded linear functional on
$L^{p'}(0, T; U')$. This implies that  there is
a $u\in \cU_p$, such that
\begin{equation}\label{5.21-eq12}
\langle \phi_T,  y_T \rangle_{Y', Y} = \int_0^T
\lan B(t)^*\phi, u(t) \rangle_{U',U} \   dt,
\q\q\forall \ \phi_T\in \wt Y
\end{equation}
and
\begin{equation}\label{8.7-eq6}
|u|_{\mathcal{U}_p}\leq C|y_T|_Y.
\end{equation}
For this control $u\in \cU_p$ and  the
corresponding  solution $y(\cdot; u)$  of
(\ref{llz51}),  by  (\ref{llz51}) and
(\ref{215}),  it is easy to show that
$$
\langle \phi_T,  y(T; u) \rangle_{Y', Y}=
\displaystyle\int_0^T \big\langle B(t)^*\phi,
u(t) \big\rangle_{U',U} dt,\q\q  \forall
\phi_T\in \wt Y,
$$
which, together with (\ref{5.21-eq12}), implies
that
$$ \langle \phi_T,  y(T; u)-y_T\rangle_{Y',
Y} = 0, \q  \forall \phi_T\in \wt Y.
$$
This deduces \eqref{00llz00}.

\smallskip

\noindent {\bf Step 3. } In this step, we prove
that for the system  (\ref{llz51}), $\cR(T;0)$
is a closed subspace of $Y$.

\smallskip

Let  $\mathbb{P}_{Y_1}$ be the projection
operator from $Y$ to  the subspace $Y_1$. Since
$Y_1$  is finite dimensional, $\mathbb{P}_{Y_1}$
is well defined and   a bounded  linear operator
on $Y$. Then for  the identity operator
$\mathbb{I}$ on $Y$, $ \mathcal{R}(T; 0)=\mathbb
P_{Y_1} \mathcal{R}(T; 0)\oplus (\mathbb
I-\mathbb P_{Y_1})\mathcal{R}(T; 0). $ Since
$\mathbb P_{Y_1} \mathcal{R}(T; 0)$ is finite
dimensional, it is closed. Furthermore, the
linear subspace $(\mathbb I-\mathbb
P_{Y_1})\mathcal{R}(T; 0)$ is also closed in
$Y$. Indeed, for any
$\{y_{T}^j\}_{j=1}^{\infty}\subseteq(\mathbb
I-\mathbb P_{Y_1})\mathcal{R}(T; 0)$
 satisfying that $\lim\limits_{j\rightarrow \infty} y_{T}^j=\widehat{y_T}\in Y$,
similar  to  (\ref{8.7-eq6}) and
(\ref{00llz00}),  there exists a  sequence  of
controls
$\{\widehat{u}_j\}_{j=1}^{\infty}\subseteq\mathcal{U}_p$,
such  that
\begin{equation}\label{llz000}
|\widehat{u}_j|_{\mathcal{U}_p}\leq
C|y_{T}^j|_Y,
\end{equation}
and for the  corresponding  mild solution
$\widehat{y}_j(\cdot)=y(\cdot; \widehat{u}_j)$
to  (\ref{llz51}),
$\widehat{y}_j(T)-y_{T}^j\in Y_1$. Therefore,
$(\mathbb I-\mathbb
P_{Y_1})\widehat{y}_j(T)=y_{T}^j$. By
(\ref{llz000}),  there   exist a  subsequence of
$\{\widehat{u}_j\}_{j=1}^{\infty}$ (still
denoted by itself) and
$\widehat{u}\in\mathcal{U}_p$,   such that as
$j\rightarrow \infty$,
\begin{eqnarray*}\left\{\!\!\!\!
\begin{array}{llll}
&\widehat{u}_j  \rightarrow
\widehat{u}\quad\quad &\mbox{
weakly in }  \mathcal{U}_p, &\mbox{  for }p\in (1, \infty);\\
\ns\ds &\widehat{u}_j \rightarrow
\widehat{u}\quad\quad &\mbox{ weakly}^*\mbox{ in
} \mathcal{U}_\infty, &\mbox{  for }p=\infty.
\end{array}
\right.
\end{eqnarray*}
Denote by $\widehat{y}$ the mild  solution to
(\ref{llz51})  associated to
$\widehat{u}\in\mathcal{U}_p$. Then it is easy
to  show  that as $j\rightarrow\infty$,
 $\widehat{y}_j(T)$  converges weakly to $\widehat{y}(T)$,
and hence, $(\mathbb I-\mathbb
P_{Y_1})\widehat{y}_j(T)$ converges weakly to
$(\mathbb I-\mathbb P_{Y_1})\widehat{y}(T)$.
Since $y_T^j$  converges strongly to
$\widehat{y_T}$ in  $Y$  and $(\mathbb I-\mathbb
P_{Y_1})\widehat{y}_j(T)=y_{T}^j$, it holds that
$ (\mathbb I-\mathbb
P_{Y_1})\widehat{y}(T)=\widehat{y_T}. $ This
implies that  $\widehat{y_T}\in (\mathbb
I-\mathbb P_{Y_1})\mathcal{R}(T; 0)$. Therefore,
$(\mathbb I-\mathbb P_{Y_1})\mathcal{R}(T; 0)$
is a closed subspace of $Y$. So is
$\mathcal{R}(T; 0)$.

\smallskip

\noindent {\bf Step 4.} In this step, we  prove
that  the codimension of the closed subspace
$\cR(T;0)$ is not greater than the dimension
$k_2$ of  $Y_1$.

\smallskip

Otherwise, there exist  linearly independent
$x_1, \cds, x_{k_2+1}\in Y$, such that for any
$\tilde x_j \in \cR(T;0)$ ($j=1, \cds, k_2+1$),
\begin{equation}\label{8.7-eq7}
\tilde x_1-x_1,\cds, \tilde x_{k_2+1}-x_{k_2+1}
\mbox{ are linearly independent. }
\end{equation}
By  (\ref{8.7-eq6}) and (\ref{00llz00}),   there
are controls $u_j\in \mathcal{U}_p$ ($j=1,\cds,
k_2+1$),  such that the corresponding mild
solutions $y_j(\cdot)=y(\cdot; u_j)$ to
(\ref{llz51}) satisfies that $ y_j(T)-x_j \in
Y_1, \mbox{ for } j=1, \cds, k_2+1. $ Meanwhile,
it follows from \eqref{8.7-eq7} that
$ y_1(T)-x_1, \cds, y_{k_2+1}(T)-x_{k_2+1}
\mbox{ are linearly independent. } $
This contradicts the fact that the dimension of
$Y_1$ is $k_2$.  Therefore, the codimension of
$\cR(T;0)$ is finite  and the assertion {\bf
(1)} holds.

\smallskip

Next,  we  prove that  ${\bf (1)}$ implies ${\bf
(2)}$.  Assume the codimension of $\cR(T;0)$ is
$k_3$  for  the system (\ref{llz51}). Then there
is a linear subspace $Y_2$ of $Y$, whose
dimension is $k_3$, such that
$ Y=\mbox{span}\big\{  Y_2,\  \cR(T;0)\big\}. $
Let $\{x_1,\cds, x_{k_3}\}$ be a basis of $Y_2$.
Then for any $x_j$ $(j=1,  \cdots, k_3)$, there
exists an $x''_j\in Y''$, such that $\langle
x''_j,  f\rangle_{Y'', Y'}= \langle  f,
x_j\rangle_{Y', Y}$, for any $ f\in Y'$. Set $
\wt Y=\bigcap\limits_{j=1}^{k_3}  \ker
\{x''_j\}. $ Then by Lemma \ref{lm1}, $\wt Y$ is
a  finite codimensional subspace of $Y'$. Also,
for any $y_T^1\in Y_2$ and $\phi_T\in \widetilde
Y$,
\begin{equation}\label{008llz}
\langle  \phi_T, y_T^1\rangle_{Y', Y}=0.
\end{equation}

Now, we prove  that    (\ref{221}) holds  for
the   above  subspace $\widetilde Y$ of $Y'$. If
(\ref{221}) fails,  then there  exists  a
sequence   $\big\{\phi_T^j
\big\}_{j=1}^{\infty}$ of  $ \widetilde  Y$,
such that  the  solution $\phi_j$ to (\ref{215})
with the final datum $\phi_j(T)=\phi_T^j$
satisfies
$$ |B(\cdot)^*\phi_j|_{L^{p'}(0, T;
U')}<\ds\frac{1}{j}|\phi^j_T|_{Y'},  \qq
\forall\ j\in\dbN.
$$ Let
$\widetilde{\phi}_T^j=\sqrt{j}\ds\frac{\phi^j_T}{|\phi_T^j|_{Y'}}$.
Then for the mild solution $\widetilde{\phi_j}$
to \eqref{215} with the final datum
$\widetilde{\phi_j}(T)=\widetilde{\phi}^j_T$, it
holds that
\begin{equation}\label{10llz}
\ds|\widetilde{\phi}_T^j|_{Y'}=\sqrt{j}\q
\mbox{and}\q
|B(\cdot)^*\widetilde{\phi_j}|_{L^{p'}(0, T;
U')}<\frac{1}{\sqrt{j}}.
\end{equation}
By \eqref{llz51} and \eqref{215}, for  any
$y_T^2\in  \mathcal{R}(T; 0)$, one can find a
control  $v(\cdot)\in  \mathcal{U}_p$,  such
that
$$
\ds\langle  \widetilde{\phi}_T^j,
y_T^2\rangle_{Y',  Y}=\int_0^T \big\langle
B(t)^*\widetilde{\phi_j}, v(t)
\big\rangle_{U',U}\    dt.
$$
This, together  with (\ref{008llz}) and
(\ref{10llz}),   implies that for any
$y_T=y_T^1+y_T^2\in Y$, $\ds\big\{\langle
\widetilde{\phi}_T^j, y_T\rangle_{Y',
Y}\big\}_{j=1}^{\infty}$ is uniformly bounded.
Hence,  $\big\{
\widetilde{\phi}_T^j\big\}_{j=1}^{\infty}$ is
uniformly bounded in  $Y'$,  but this
contradicts  (\ref{10llz}). Hence,  (\ref{221})
holds for any $\phi_T \in \wt Y$.  \endpf

\medskip

\medskip

In general,  it is hard to find the finite
codimensional  subspace $\widetilde{Y}$  of $Y'$
in the assertion {\bf (2)} of Theorem
\ref{t2.1}. Hence, we give another equivalent
criterion  for the finite codimensional exact
controllability, where {\it a priori} estimate
holds on the whole space $Y'$.
\begin{theorem}\label{t3}
The following two assertions are equivalent:

\medskip

\noindent  {\bf (1) } There is a  compact
operator $G$ from  $Y'$ to a Banach space $X$,
such that   any  solution $\phi$ to
$(\ref{215})$ satisfies
\begin{equation}\label{5.21-eq1}
|\phi_T|_{Y'}\leq
C\big(|B(\cdot)^*\phi|_{L^{p'}(0, T; U')} +
|G\phi_T|_{X}\big),  \quad\quad\forall\
\phi_T\in Y'.
\end{equation}

\noindent  {\bf (2) } There is a finite
codimensional subspace $\wt Y$ of   $Y'$,  such
that  any  solution $\phi$ to $(\ref{215})$
satisfies
\begin{equation}\label{22111}
|\phi_T|_{Y'}\leq C |B(\cdot)^*\phi|_{L^{p'}(0,
T; U')},\quad\quad\forall\ \phi_T\in \wt Y.
\end{equation}
\end{theorem}

\noindent {\bf  Proof.} First,  we prove that
{\bf (1)} implies  {\bf (2)}.  The proof is
divided into three parts.

\smallskip

\noindent {\bf Step 1.} In this step, we  prove
that the following subspace  $\cE$ of $Y'$  is
finite dimensional:
$$
\cE=\Big\{ \phi_T\in Y'\ \Big|\  \mbox{the
corresponding mild solution $\phi$ to
(\ref{215})  satisfies that } B(\cdot)^*\phi=0
\Big\}.
$$
 Indeed, let
$\{\phi^j_T\}_{j=1}^\infty\subseteq\cE$ with
$|\phi^j_T|_{Y'}=1$ for every  $j\in\dbN$. Then
there exist  a $\tilde\phi_T \in Y'$ and
subsequence of $\{\phi^j_T\}_{j=1}^\infty$
(still denoted by itself),  such that
$$ \phi^j_T\rightarrow  \tilde\phi_T\q\q  \mbox{
weakly}^*\mbox{ in }
  Y' ,\quad\quad \mbox{ as }j\rightarrow
+\infty.$$ Hence, $
\lim\limits_{j\to\infty}G\phi^{j}_T=G\tilde\phi_T
\mbox{ in } X. $ This, together with
\eqref{5.21-eq1}, implies that
$\{\phi_T^j\}_{j=1}^\infty$ is strongly
convergent
 in $Y'$ and therefore, $\cE$ is a
finite dimensional space.

\smallskip

\noindent {\bf Step  2.} We find a suitable
$\hat{\phi}_T\in \cE $  with $\hat{\phi}_T\neq
0$, by assuming the following  (\ref{**LLLL})
fails.

\smallskip

Denote by $\mathbb P_\cE$ the projection
operator from $Y'$ to the subspace $\cE$. Since
$\cE$ is finite dimensional, $\mathbb P_\cE$ is
well defined and a bounded linear operator.
  In the following,   we
prove that \begin{equation}\label{**LLLL}
|\phi_T|_{Y'}\leq C |B(\cdot)^*\phi|_{L^{p'}(0,
T; U')}, \mbox{  for all } \phi_T\in
(\mathbb{I}-\mathbb P_\cE)Y',
\end{equation}where $\mathbb{I}$ denotes the identity
operator  on $Y'$ and $\phi$ is the mild
solution to (\ref{215}) with the terminal  value
$\phi_T$. Otherwise,
 there exists a  sequence
$\{\phi_{T,  j}\}_{j=1}^\infty$ of $(\mathbb
I-\mathbb P_\cE)Y'$ with $|\phi_{T,  j}|_{Y'}=1$
for any $j\in\dbN$, such that  the solution
$\phi_j$ to  (\ref{215}) with the final datum
$\phi_j(T)=\phi_{T, j}$ satisfies
\begin{equation}\label{5.21-eq4}
|B(\cdot)^*\phi_j|_{L^{p'}(0, T;
U')}<\frac{1}{j}.
\end{equation}
Then there exist a  $\hat{\phi}_T\in Y'$ and
subsequence
 of
$\{\phi_{T,  j}\}_{j=1}^\infty$ (still denoted
by itself),  such that
%\begin{equation}\label{8.7-eq12}
$$\ds\phi_{T, j} \rightarrow  \hat{\phi}_T\q\q  \mbox{
weakly}^*\mbox{ in }
  Y' ,\quad\quad \mbox{ as }j\rightarrow
+\infty.
$$
%\end{equation}
%
By (\ref{5.21-eq4}), one has that
$B(\cdot)^*\hat{\phi}=0$, where $\hat{\phi}$ is
the mild solution to (\ref{215}) with the final
datum  $\hat\phi(T)=\hat{\phi}_T$. This  implies
that $\hat{\phi}_T\in\cE$.
Also,\begin{equation}\label{5.21-eq6}
\lim_{j\to\infty}G\phi_{T, j}=G\hat{\phi}_T\qq
\mbox{ in } X.
\end{equation}
It follows from (\ref{5.21-eq1}) that
$ |\phi_{T,  j}|_{Y'}\leq
C\big(|B(\cdot)^*\phi_j|_{L^{p'}(0, T; U')} +
|G\phi_{T, j}|_{X}\big). $
Hence, by (\ref{5.21-eq4}),  for any $j>2C$, one
has that $ |G\phi_{T,  j}|_{X}\geq\ds
\frac{1}{2C},$
which,  together with   (\ref{5.21-eq6}),
indicates that  $\hat{\phi}_T\neq  0$.

\smallskip

\noindent {\bf Step  3.} In this step,  we prove
that $\hat{\phi}_T\in (\mathbb I-\mathbb
P_\cE)Y'$.

\smallskip

By (\ref{5.21-eq1}), for the above $\{\phi_{T,
j}\}_{j=1}^\infty\subseteq (\mathbb I-\mathbb
P_\cE)Y'$,
\begin{eqnarray*}
\begin{array}{ll}\ds
|\phi_{T,  n}-\phi_{T,  m}|_{Y'}\\
\ns\ds\leq
C\Big(|B(\cdot)^*\phi_{n}-B(\cdot)^*\phi_{m}|_{L^{p'}(0,
T; U')} + |G\phi_{T, n}-G\phi_{T, m}|_{X}\Big)\\
\ns \ds \leq C\(\frac{1}{n} + \frac{1}{m} +
|G\phi_{T,  n}-G\phi_{T,   m}|_{X}\),\q \mbox{
for  any  }m, n\in\dbN,
\end{array}
\end{eqnarray*}
which implies that $\{\phi_{T,
j}\}_{j=1}^\infty$ is a Cauchy sequence in $Y'$.
Since $(\mathbb I-\mathbb P_\cE)Y'$ is closed,
$\hat{\phi}_T\in (\mathbb I-\mathbb P_\cE)Y'$.
This  contradicts the fact that $\hat{\phi}_T\in
\cE$ and $\hat{\phi}_T\neq 0$. Hence,
(\ref{22111}) holds,    provided that  $\phi_T$
belongs to the finite codimensional  subspace
$\wt Y=(\mathbb{I}-\mathbb{P}_{\cE})Y'.$

\smallskip

Next,   we  prove that {\bf (2)} implies {\bf
(1)}. Assume that there is a finite dimensional
subspace $Y_3$  of $Y'$,  such that $
Y'=\mbox{span}\{ \widetilde{Y},  Y_3 \}. $
Denote by $\mathbb P_{Y_3}$  and $\mathbb
P_{\widetilde{Y}}$ the projections from $Y'$ to
$Y_3$ and $\widetilde{Y}$, respectively. Also,
for any $\phi_T\in Y'$,  denote  by
$\mathcal{F}(\phi_T)$  the associated solution
$\phi$ to (\ref{215}).  Then by (\ref{22111}),
for any $\phi_T\in Y'$,  it holds that
\begin{eqnarray}\label{12llz}
\begin{array}{ll}
\ds|\phi_T|_{Y'}\leq |\mathbb P_{Y_3}
\phi_T|_{Y'}+|\mathbb P_{\widetilde
Y}\phi_T|_{Y'}
\leq |\mathbb P_{Y_3} \phi_T|_{Y'}+C\big|B(\cdot)^*\mathcal{F}(\mathbb P_{\widetilde{Y}}\phi_T)\big|_{L^{p'}(0, T; U')}&\\
\ns\ds\leq C\big( |\mathbb P_{Y_3}
\phi_T|_{Y'}+\big|B(\cdot)^*\mathcal{F}(\mathbb
P_{Y_3}\phi_T)\big|_{L^{p'}(0, T; U')}\big)+
C\big|B(\cdot)^*\phi\big|_{L^{p'}(0, T; U')}.&
\end{array}
\end{eqnarray}
Define a linear operator:
$$G: Y'\rightarrow
Y'\times L^{p'}(0, T; U'),\q  \q
G(\phi_T)=\big(\mathbb P_{Y_3}\phi_T,
B(\cdot)^*\mathcal{F}(\mathbb
P_{Y_3}\phi_T)\big), \q \forall\  \phi_T\in Y'.
$$
Since  $\mathbb P_{Y_3}$ is compact,   $G$    is
also  compact from $Y'$  to  the Banach space
$X=Y'\times L^{p'}(0, T; U')$ and  therefore,
(\ref{5.21-eq1})  follows from (\ref{12llz}).
\endpf

\section{Finite codimensional approximate
controllability}\label{sec-ap}

 This
section is devoted to a characterization of the
finite codimensional approximate
controllability.

\medskip

The  main result of this section is stated as
follows.
\begin{theorem}\label{t1}
Suppose that ${\bf (A)}$ holds.  Then the
following three assertions are equivalent:

\medskip

\noindent  {\bf (1)}  The system $(\ref{llz51})$
is finite  codimensional approximately
controllable in $Y$.

\medskip

\noindent {\bf  (2)} \mbox{span}M $($in
$(\ref{llz53}))$  is a finite codimensional
subspace of $Y$.

\medskip

\noindent  {\bf (3)} There is  a finite
dimensional subspace $\widehat Y$ of $Y'$, such
that for any solution $\phi$ to $(\ref{215})$,
$$
B(\cdot)^*\phi=0,\ \   \mbox{  if and only if
}\quad \phi_T\in \widehat  Y. $$

\end{theorem}

\noindent {\bf Proof.}   First, by Definition
\ref{d1}  and (\ref{llz53}),  it is obvious that
{\bf (1)} and {\bf (2)} are equivalent.

\smallskip

Next,  we prove that {\bf (2)} implies  {\bf
(3)}. Define a linear operator $\mathbb{L}:
\mathcal{U}_p\rightarrow Y$ as $
\mathbb{L}(u(\cdot))=y(T; u), \    \forall
u(\cdot)\in\mathcal{U}_p, $ where $y(\cdot; u)$
is the mild solution to (\ref{llz51}) associated
to $u$.  Denote by $\mathcal R(\mathbb{L})$ the
range of  the operator $\mathbb{L}.$  Then
$\mathbb{L}$ is a bounded linear operator and
span$M=\overline{\mathcal  R(\mathbb{L})}$.
Therefore,   {\bf (2)}  means that
 $\overline{\mathcal R(\mathbb{L})}$ is a finite codimensional subspace of $Y$.
Also, it is easy to show that the adjoint
operator $\mathbb{L}^*:  Y'\rightarrow
(\mathcal{U}_{p})'$ of $\mathbb{L}$ is  $
\mathbb{L}^*(\phi_T)= B(\cdot)^*\phi, $ where
$\phi$ is the mild solution to (\ref{215})
associated to $\phi_T\in Y'$. Since
$$\mbox{ker}(\mathbb{L}^*)=
\big(\overline{\mathcal
R(\mathbb{L})}\big)^\perp\=\big\{ g\in Y'\
\big|\ \langle   g, x\rangle_{Y', Y}=0,\
\forall\ x\in \overline{\mathcal R(\mathbb{L})}
\big\},
$$
it suffices to show that
$\big(\overline{\mathcal
R(\mathbb{L})}\big)^\perp$ is finite
dimensional.

\smallskip

 Since $\overline{\mathcal
R(\mathbb{L})}$ is finite codimensional in $Y$,
there exists a finite dimensional subspace $Y_0$
of  $ Y$, such that $Y_0\oplus
\overline{\mathcal R(\mathbb{L})}=Y$.  Denote by
$m_0$ and $\dbP_{Y_0}$
 the dimension of $Y_0$  and the projection
operator from $Y$ to $Y_0$,  respectively. Since
$Y_0$ is finite dimensional, $\mathbb{P}_{Y_0}$
is well defined and a bounded linear operator.
If $\big(\overline{\mathcal
R(\mathbb{L})}\big)^\perp$ is infinite
dimensional in $Y'$,    there exist linearly
independent $f_1, f_2, \cdots,
f_{m_0+1}\in\big(\overline{\mathcal
R(\mathbb{L})}\big)^\perp$.  For the above given
$f_j$ ($j=1, \cdots, m_0+1$), define  a bounded
linear functional $\tilde{f}_{j}$ on $Y_0$  as
the limitation  of $f_j$ on $Y_0$.  Since
$\tilde{f}_{j}\in Y_0'$  $(j=1, \cdots,
m_{0}+1)$,
 the dimension of
$\span\{\tilde{f}_{1}, \cdots,
\tilde{f}_{m_0+1}\}$ is not larger than $m_0$.
Hence, there exists  a nonzero vector
$(a_1,\cds,a_{m_0+1})^\top\in\dbR^{m_0+1}$,
such that $
\sum\limits_{j=1}^{m_0+1}a_j\tilde{f}_{j}(x)=0,\
\forall\ x\in Y_0. $ Notice that   for  any
$x\in Y$,
$$
\sum_{j=1}^{m_0+1}a_jf_{j}(x)=\sum_{j=1}^{m_0+1}a_jf_{j}\big(\mathbb{P}_{Y_0}x+(\mathbb{I}-\mathbb{P}_{Y_0})x\big)=\sum_{j=1}^{m_0+1}a_jf_{j}(\mathbb{P}_{Y_0}x)
=\sum_{j=1}^{m_0+1}a_j\tilde{f}_{j}(\mathbb{P}_{Y_0}x)=0.
$$
This implies that
$\sum\limits_{j=1}^{m_0+1}a_jf_{j}=0$ in $Y'$.
This contradicts the linear independence   of
$f_1, \cdots, f_{m_0+1}$.

Finally,  we prove that {\bf (3)} implies  {\bf
(2)}. By {\bf (3)},
$\mbox{ker}(\mathbb{L}^*)=\big(\overline{\mathcal
R(\mathbb{L})}\big)^\perp$ is finite dimensional
in $Y'$.  Assume that its dimension  is $m_1$.
If $\overline{\mathcal R(\mathbb{L})}$ is not
finite codimensional, then there exist
 $x_j\notin
\overline{\mathcal R(\mathbb{L})}$
$(j=1,\cds,m_1+1)$, which are linearly
independent  in $Y$. Let $Y_1=\span\{x_1,
\cdots, x_{m_1+1}\}\oplus \overline{\mathcal
R(\mathbb{L})}$.  Then $Y_1$ is closed and for
any $x\in Y_1$,  there  exists  a vector
$(a_{1, x},  \cdots, a_{m_1+1,
x})^\top\in\dbR^{m_1+1}$,  such that $
x=\sum\limits_{j=1}^{m_1+1}a_{j,x}x_j + \tilde
x$ for some $\tilde x\in \overline{\mathcal
R(\mathbb{L})}$. Define a functional $g_j$
$(j=1, \cdots, m_1+1)$ on $Y_1$ as follows:
\begin{center}
$ g_j(x) =a_{j,x},\qq \forall\
x=\sum\limits_{j=1}^{m_1+1}a_{j,x}x_j + \tilde
x\in Y_1. $
\end{center}
It is obvious that $g_j$ is linear and bounded.
By the Hahn-Banach theorem, $g_j$ has an
extension $\tilde g_j\in Y'$,  such that
$g_j(x)=\tilde g_j(x)$ for all $x\in
\overline{\mathcal R(\mathbb{L})}$. Hence,
$\{\tilde g_j\}_{j=1}^{m_1+1}\subseteq
\big(\overline{\mathcal
R(\mathbb{L})}\big)^\perp$ and $\tilde g_1,
\cdots, \tilde g_{m_1+1}$ are linearly
independent in $Y'$. This contradicts the fact
that the dimension of $\big(\overline{\mathcal
R(\mathbb{L})}\big)^\perp$ is $m_1$.
\endpf

\medskip

\medskip

Furthermore,   suppose that  the equation
(\ref{215}) satisfies the following forward
uniqueness:
$$
{\bf (U)}\q \mbox{   For any solution } \phi
\mbox{ to }  (\ref{215}), \mbox{ if } \phi(0)=0,
\mbox{ then } \phi(\cdot)=0 \mbox{ in } [0, T].
$$
Then,  it  is easy to show  the following
equivalence result.
\begin{corollary}
Assume that ${\bf (U)}$ holds. Then the
following two assertions are equivalent:

\smallskip

\noindent  {\bf  (1) }  There is  a finite
dimensional subspace $\widehat Y$ of $Y'$, such
that for any solution $\phi$ to $(\ref{215})$,
$$
B(\cdot)^*\phi=0,\ \   \mbox{  if and only if
}\quad \phi_T\in \widehat  Y. $$

\smallskip

\noindent  {\bf  (2) }  There is  a finite
dimensional subspace $\widehat Y$ of $Y'$, such
that for any solution $\phi$ to \eqref{215},
$$
B(\cdot)^*\phi=0,\ \   \mbox{  if  and  only if
}\quad \phi(0)\in \widehat Y. $$

\end{corollary}

\section{Applications to  optimal control problems with state constraints }\label{sec-op}

In this section,  as applications of results on
the finite codimensional  exact  controllability
in Section \ref{sec-ex},    we  study the finite
codimensionality  {\bf  (H)} for the optimal
control problem {\bf  (P)}  with state
constraints in Section \ref{sss3}.

\medskip

In the following,  suppose that
$(\overline{u}(\cdot), \overline{y}(\cdot))$ is
an optimal
 pair of   the  optimal control  problem  ${\bf (P)}$.  We will study its necessary conditions    in two different cases.

 \medskip

\noindent {\bf  Case 1. }  The optimal control
problem {\bf (P)}    {\bf
 without control constraints}.

 \medskip

First,   we give the following hypotheses:

\medskip

\noindent ${\bf (A_{11})}$  Let  $f: [0,
T]\times Y\times U\rightarrow Y$ and $f^0: [0,
T]\times Y\times U\rightarrow \dbR$ be strongly
measurable with respect to $t$ in $(0, T)$, and
continuously Fr\'echet differentiable with
respect to  $(y, u)$ in $Y\times U$ with $f(t,
\cdot, \cdot), f_y(t, \cdot, \cdot),  f_u(t,
\cdot, \cdot), f^0(t, \cdot, \cdot), f^0_y(t,
\cdot, \cdot)$ and $f_u^0(t, \cdot, \cdot)$
being continuous. Moreover,   for any
$u(\cdot)\in\mathcal{U}_2$ and $y(\cdot)\in
C([0, T]; Y)$,
$$
\begin{array}{ll}
f_y(\cdot, y(\cdot), u(\cdot))\in L^1(0,  T;
\mathcal{L}(Y)),\ \ f^0_y(\cdot, y(\cdot),
u(\cdot))\in L^1(0,  T; Y'),\  \ f^0(\cdot,
y(\cdot),u(\cdot))\in L^1(0, T), \\\ns\ds
f_u(\cdot, y(\cdot), u(\cdot))\in L^2(0,  T;
\mathcal{L}(U; Y))\  \mbox{  and  }\
f^0_u(\cdot, y(\cdot), u(\cdot))\in L^2(0,  T;
U').
\end{array}
$$

\noindent  ${\bf (A_{12})}$  $p=2$  and
$\widetilde{U}=U$.

\medskip

\noindent Under the above assumptions,
$\mathcal{U}[0, T]=\mathcal{U}_2=L^2(0, T; U)$.
For any $T>0$, $y_0\in Y$ and
$u(\cdot)\in\mathcal{U}_2$, $(\ref{61})$ admits
a
 mild solution $y(\cdot)=y(\cdot; y_0, u)\in C([0, T]; Y)$
 corresponding  to the  initial value $y_0$ and the control $u$. Also, we  write
 $J(y_0, u(\cdot))=J(u(\cdot), y(\cdot))$.

\medskip

Consider the  following linear control system:
\begin{eqnarray}\label{64}
\left\{
\begin{array}{ll}
\ds\widetilde\xi_t(t)=A\wt\xi(t)+f_y(t,
\overline{y}(t), \overline{u}(t))\wt\xi(t)+
f_u(t, \overline{y}(t), \overline{u}(t))v(t),  & t\in (0, T],\\
\ns\ds\wt\xi(0)=y^0,
\end{array}
\right.
\end{eqnarray}
where $v(\cdot)\in  \mathcal{U}_2$ and set
$$
\begin{array}{ll}
\ds M_1=\Big\{ \wt\xi(T)-y^1\in  Y\ \Big| \
\wt\xi\mbox{ is the solution to }
(\ref{64})\mbox{ with some
}v(\cdot)\in\mathcal{U}_2 \mbox{ satisfying that
} |v|_{\mathcal{U}_2}\leq 1 \\\ns\ds
\q\q\q\q\q\q\q\q\q\q\  \mbox{ and } (y^0,y^1)\in
S \Big\}.
\end{array}
$$
Further,  let $\psi_1$ be the mild solution to
the following  equation:
\begin{eqnarray}\label{65}
\left\{\3n
\begin{array}{ll}
\ds\psi_{1, t}(t)=-A^*\psi_1(t)-\psi^0_1
f^0_y(t, \overline{y}(t), \overline{u}(t))-
f_y(t, \overline{y}(t), \overline{u}(t))^*\psi_1(t),  & t\in (0, T],\\
\ns\ds\psi_1(T)=\psi^1_T,
\end{array}\right.
\end{eqnarray}
where $(\psi^0_1, \psi^1_T)\in \dbR\times Y'$.
Then,  similar to  \cite{2}, using a convex
variation  technique, we can derive  a necessary
condition for the optimal pair
$(\overline{u}(\cdot), \overline{y}(\cdot))$  as
follows.
\begin{proposition} \label{tt1}
Assume that  ${\bf (A_{11})}$ and ${\bf
(A_{12})}$  hold,  and  $M_1$  is {\bf finite
codimensional} in $Y$.  Then there exists a pair
$(\psi^0_1, \psi^1_T)\in \dbR\times Y'$, such
that for the corresponding  solution $\psi_1$ to
$(\ref{65})$,  $(\psi^0_1, \psi_1(\cdot))\neq
(0, 0)$ and
\begin{equation}\label{011llz}
\ds f_u(t, \overline{y}(t),
\overline{u}(t))^*\psi_1(t)+\psi^0_1 f^0_u(t,
\overline{y}(t),
\overline{u}(t))=0,\quad\mbox{a.e. }t\in (0, T).
\end{equation}\end{proposition}
\noindent  {\bf  Proof. }  The proof is divided
into   four steps.

\smallskip

\noindent {\bf Step 1.  } For any $\e>0$,
$y_0\in Y$ and $u(\cdot)\in \mathcal{U}_2$, set
$\hat J(y_0, u(\cdot)) =J(y_0,
u(\cdot))-J(\overline{y}_0,
\overline{u}(\cdot))$  and
$$
J_\e(y_0, u(\cdot))= \Big\{\big[d_S(y_0, y(T;
y_0, u))\big]^2 +\big[\hat J(y_0,
u(\cdot))+\e\big]^2\Big\}^{1/2},
$$
where $\overline{y}_0=\overline{y}(0)$ and
$d_S(y_0, y_1) =\inf\limits_{(y^0, y^1)\in
S}[|y_0-y^0|^2+|y_1-y^1|^2]^{1/2}$. Then
$J_\e(\cdot, \cdot)$ is continuous on
$Y\times\mathcal{U}_2$ and $$
J_\e(\overline{y}_0, \overline{u}(\cdot))=\e\leq
\inf\limits_{(y_0, u(\cdot))\in
Y\times\mathcal{U}_2} J_\e (y_0, u(\cdot))+\e.
$$ By the Ekeland variational principle,  there
exists  a pair $(y_0^\e,  u_\e(\cdot))\in
Y\times\mathcal{U}_2$, such that
 $$
 J_\e(y_0^\e, u_\e(\cdot))\leq J_\e(\overline{y}_0, \overline{u}(\cdot)),\quad\quad
 |y_0^\e-\overline{y}_0|_Y+|u_\e(\cdot)-\overline{u}(\cdot)|_{\mathcal{U}_2}
 \leq \sqrt\e,
 $$
 and
 \begin{equation}\label{66}
 -\sqrt{\e}[
 |y_0^\e-y_0|_Y+
 |u_\e(\cdot)-v(\cdot)|_{\mathcal{U}_2}]
 \leq J_\e(y_0, v(\cdot))
 -J_\e(y_0^\e, u_\e(\cdot)),\quad \forall\
 (y_0, v(\cdot))\in Y\times\mathcal{U}_2.
 \end{equation}

 \smallskip

\noindent {\bf Step 2.} Set
$y_\e(\cdot)=y(\cdot; y_0^\e, u_\e)$  and for
any $\rho>0$, $\nu\in Y$ and $v(\cdot)\in
\mathcal{U}_2$, write $
u_\e^\rho(\cdot)=u_\e(\cdot)+\rho v(\cdot), \
y_0^{\e, \rho}=y_0^\e+\rho\nu,\
  y_\e^\rho(\cdot)=y(\cdot;  y_0^{\e, \rho},  u_\e^\rho), $
$$ z_\e(\cdot)=\lim\limits_{\rho\rightarrow 0}\displaystyle\frac{y^\rho_\e(\cdot)-y_\e(\cdot)}{\rho}\q
  \mbox{ and }  \q
 z^0_\e=\lim\limits_{\rho\rightarrow 0}\displaystyle\frac{\hat J(y_0^{\e, \rho}, u^\rho_\e(\cdot))-\hat J(y_0^\e, u_\e(\cdot))}{\rho}.
 $$
Then it is easy to check that $z_\e(\cdot)$ and
$z^0_\e$, respectively,   satisfy that
\begin{eqnarray*}
\left\{
\begin{array}{ll}
\ds z_{\e, t}(t)=Az_\e(t)+f_y(t, y_\e(t),
u_\e(t))z_\e(t)+
f_u(t, y_\e(t), u_\e(t)) v(t),  &\quad t\in (0, T],\\
\ns\ds z_\e(0)=\nu,
\end{array}\right.
\end{eqnarray*}
and
$$
z^0_\e= \int^T_0 \Big[\langle f^0_y(t, y_\e(t),
u_\e(t)), z_\e(t)\rangle_{Y', Y} +\langle
f^0_u(t, y_{\e}(t), u_\e(t)), v(t)\rangle_{U',
U}\Big]dt.
$$
Furthermore, by the definition of $J_\e(\cdot,
\cdot)$,  as $\rho\rightarrow 0$,
\begin{eqnarray}\label{61llz}
\begin{array}{rl}
 &\displaystyle\frac{J_\e(y_0^{\e, \rho},   u_\e^\rho(\cdot))-J_\e(y_0^\e, u_\e(\cdot))}{\rho}
 \rightarrow
 z^0_\e\psi^0_{1, \e}+\langle \psi^1_{1, \e},
 \nu\rangle_{Y', Y}+\langle \psi^2_{1, \e},
 z_\e(T)\rangle_{Y', Y},
 \end{array}
 \end{eqnarray}
where
$$\psi^0_{1, \e}=
\displaystyle\frac{\hat J(y_0^\e,
u_\e(\cdot))+\e}{J_\e(y_0^\e, u_\e(\cdot))}, \
\  \psi^1_{1, \e}=\frac{d_S(y_0^\e, y_\e(T))
a_\e}{J_\e(y_0^\e, u_\e(\cdot))} \    \mbox{ and
}\   \psi^2_{1, \e}=\frac{d_S(y_0^\e, y_\e(T))
b_\e}{J_\e(y_0^\e, u_\e(\cdot))}$$ with
$|a_\e|_{Y'}^2+|b_\e|^2_{Y'}=1$, $|\psi_{1,
\e}^0|^2+|\psi^1_{1, \e}|_{Y'}^2+|\psi^2_{1,
\e}|_{Y'}^2=1,$ and
\begin{equation}\label{L11}
\langle a_\e, y^0-y_0^\e\rangle_{Y', Y} +\langle
b_\e, y^1-y_\e(T)\rangle_{Y', Y}\leq  0,  \mbox{
for any }(y^0, y^1)\in S.
\end{equation}

On the other hand, by (\ref{66}),
\begin{equation}\label{610}
J_\e(y_0^{\e, \rho}, u_\e^\rho(\cdot))-
J_\e(y_0^\e, u_\e(\cdot))\geq -\sqrt{\e}\rho
\big(|\nu|_Y+|v(\cdot)|_{\mathcal{U}_2}\big).
\end{equation}
(\ref{61llz}) and (\ref{610}) imply that for any
$\e>0$,
\begin{equation}\label{611}
-\sqrt{\e}
\big[|\nu|_Y+|v(\cdot)|_{\mathcal{U}_2}\big]
\leq z^0_\e\psi^0_{1, \e}+\langle \psi^1_{1,
\e}, \nu\rangle_{Y', Y}+\langle \psi^2_{1, \e},
z_\e(T)\rangle_{Y', Y}.
\end{equation}

\smallskip

\noindent {\bf Step 3. }  Without loss  of
generality,  we  assume that  as $\e\rightarrow
0$, $\psi_{1, \e}^1 \rightarrow \psi_1^1 \mbox{
weakly}^*\mbox{ in }Y',\ \psi_{1, \e}^2
\rightarrow  \psi_T^1 \mbox{ weakly}^*\mbox{  in
}Y',\ \psi^0_{1,  \e}\rightarrow \psi^0_1\mbox{
in }\dbR, $    $z^0_\e\rightarrow z^0$ in
$\dbR$,   and  $\sup\limits_{t\in [0,
T]}|z_\e(t)- z(t)|_Y\rightarrow 0$. Then by
(\ref{611}),
\begin{equation}\label{00LLZ0}
\ds z^0 \psi^0_{1}+\langle \psi^1_{1},
 \nu\rangle_{Y', Y}+\langle \psi^1_{T},
 z(T)\rangle_{Y', Y}\geq 0, \mbox{ for any }
(\nu, v(\cdot))\in Y\times\mathcal{U}_2,
\end{equation}
where  $z^0$ and $z(\cdot)$ satisfy,
respectively,  that
$$
\ds z^0= \int^T_0 \Big[\langle f^0_y(t,
\overline{y}(t), \overline{u}(t)),
z(t)\rangle_{Y', Y} +\langle f^0_u(t,
\overline{y}(t), \overline{u}(t)),
v(t)\rangle_{U', U}\Big]dt
$$
and
\begin{eqnarray}\label{67}
\left\{
\begin{array}{ll}
\ds z_{t}(t)=Az(t)+f_y(t, \overline{y}(t),
\overline{u}(t))z(t)+
f_u(t, \overline{y}(t), \overline{u}(t)) v(t),  &\quad t\in (0, T],\\
\ns\ds z(0)=\nu.
\end{array}\right.
\end{eqnarray}

On the other hand,  let $\psi_1$ be the solution
to (\ref{65})  associated  to the above
$(\psi^0_1,  \psi_T^1)\in\dbR\times Y'$.  By
(\ref{65}), (\ref{00LLZ0}) and (\ref{67}),
choose $\nu=0$ and it is  easy to find that
$$
\ds\int_0^T \Big[\langle\psi^0_1 f^0_u(t,
\overline{y}(t), \overline{u}(t)),
v(t)\rangle_{Y', Y}+ \langle f_u(t,
\overline{y}(t), \overline{u}(t))^*\psi_1,
v(t)\rangle_{Y', Y}\Big] dt=0,\ \forall\
v(\cdot)\in \mathcal{U}_2.
$$
This implies  the necessary  condition
(\ref{011llz}).

\smallskip

\noindent {\bf Step 4. }   The finite
codimensionality of $M_1$ is given to guarantee
that $(\psi^0_1, \psi_1(\cdot))\neq (0, 0)$.
Indeed,
 if $\psi^0_1=0$, there exists a $\delta_0>0$,  such  that for sufficiently small $\e>0$,
$ |\psi^1_{1, \e}|^2_{Y'}+|\psi^2_{1,
\e}|_{Y'}^2\geq \delta_0. $ Also, by
(\ref{611}) and (\ref{L11}), it follows that for
any $(\nu, v(\cdot))\in Y \times \mathcal{U}_2$
with $|\nu|_Y\leq 1$ and
$|v(\cdot)|_{\mathcal{U}_2}\leq 1$,
\begin{eqnarray*}
\begin{array}{ll}
\ds\langle  \psi_{1, \e}^1,
\nu-y^0+\overline{y}_0\rangle_{Y', Y}
+\langle\psi^2_{1, \e},
z(T)-y^1+\overline{y}(T)\rangle_{Y', Y}&\\\ns\ds
\geq
-\sqrt{\e}\big[|\nu|_Y+|v(\cdot)|_{\mathcal{U}_2}\big]
+\psi^0_{1, \e} z^0+\psi^0_{1, \e}(z^0-z^0_\e)
-\langle \psi^1_{1,\e},
y^0-\overline{y}_0\rangle_{Y', Y} &\\\ns\ds
\quad -\langle \psi^2_{1, \e},
y^1-\overline{y}(T)\rangle_{Y', Y} +\langle
\psi^2_{1, \e},  z(T)-z_\e(T)\rangle_{Y',
Y}&\\\ns\ds \geq -2\sqrt{\e}-|\psi^0_{1, \e}
z^0|-|z^0-z^0_\e|-
|y_0^\e-\overline{y}_0|_Y-|y_\e(T)-\overline{y}(T)|_Y
-|z_\e(T)-z(T)|_Y&\\\ns\ds
\=-\delta_\e\rightarrow 0,  \mbox{  as
}\e\rightarrow 0.
\end{array}
\end{eqnarray*}
Then,  by  Proposition 3.5 and Lemma 3.6 in
\cite{2},   if   $M_1$ is  finite  codimensional
in $Y$, $(\psi_1^0, \psi_1^1,  \psi^1_T)\neq (0,
0,  0)$. This implies that $(\psi^0_1,
\psi_1(\cdot))$ is  nontrivial, since
$\psi_1^1=-\psi_1(0)$.    \endpf

\medskip

In the following,  we study the finite
codimensionality  of  the set  $M_1$  under
   fixed endpoint
constraints. Set $S=\big\{(y^0, y^1)\big\}$,
where  $y^0,  y^1\in Y$ are arbitrarily  given.
Then   the solution $\wt\xi$ to  \eqref{64}  can
be rewritten  as $\wt\xi=z_1+\eta_1$, where
$z_1$  and  $\eta_1$ solve, respectively, that
\begin{eqnarray}\label{618}
\left\{
\begin{array}{ll}
\ds  z_{1,t}(t)=Az_1(t)+f_y(t, \overline{y}(t),
\overline{u}(t))z_1(t)+
f_u(t, \overline{y}(t), \overline{u}(t))v(t),  &\quad t\in (0, T],\\
\ns\ds z_1(0)=0,
\end{array}
\right.
\end{eqnarray}
and
\begin{eqnarray}\label{619}
\left\{
\begin{array}{ll}
\ds\eta_{1, t}(t)=A\eta_1(t)+f_y(t, \overline{y}(t), \overline{u}(t))\eta_1(t),  &\quad t\in (0, T],\\
\ns\ds\eta_1(0)=y^0.
\end{array}
\right.
\end{eqnarray}
Set
\begin{equation}\label{llz*}
M_2=\Big\{  z_1(T)\in Y\ \Big|\ z_1\mbox{ is the
solution to } (\ref{618})\mbox{ with some
}v(\cdot)\in\mathcal{U}_2 \mbox{ satisfying  }
|v|_{\mathcal{U}_2} \leq 1 \Big\}.
\end{equation}
Then $M_1=M_2+\big\{\eta_1(T)-y^1\big\}$  and
therefore,    $M_1$
  is finite codimensional in $Y$, if and only if
 $M_2$ is  finite  codimensional in $Y$.

\medskip

By Theorems \ref{t2}-\ref{t3},  one has the
following equivalent   assertions   on the
finite  codimensionality  of $M_2$.

\begin{corollary}\label{LLLL}  The following  assertions are equivalent:

\medskip

\noindent  {\bf (1)}  The system $(\ref{618})$
is finite  codimensional exactly controllable in
$Y$.

\medskip

\noindent  {\bf (2)} There is a finite
codimensional subspace $\wt Y\subseteq Y'$, such
that  any  solution $\phi$ to the equation:
\begin{eqnarray}\label{002llz}
\left\{
\begin{array}{ll}
\ds \phi_t(t)=-A^*\phi(t)-f_y(t,
\overline{y}(t), \overline{u}(t))^*\phi(t),
&\quad t\in (0, T],\\
\ns\ds\phi(T)=\phi_T,
\end{array}\right.
\end{eqnarray}
satisfies that
$$
|\phi_T|_{Y'}\leq C |f_u(\cdot,
\overline{y}(\cdot),
\overline{u}(\cdot))^*\phi|_{L^{2}(0, T;
U')},\quad\quad\forall \phi_T\in \wt Y.
$$

\medskip

\noindent {\bf (3)}  There  is  a  compact
operator $G$  from  $Y'$ to a Banach space $X$,
such  that  any  solution $\phi$ to
$(\ref{002llz})$  satisfies that
\begin{equation}\label{20llz}
|\phi_T|_{Y'}\leq C\big(|f_u(\cdot,
\overline{y}(\cdot),
\overline{u}(\cdot))^*\phi|_{L^{2}(0, T; U')} +
|G\phi_T|_{X}\big),\quad\quad\forall \phi_T\in
Y'.
\end{equation}

\noindent{\bf (4)} The set  $M_2$ $($defined in
$(\ref{llz*}))$ is  finite codimensional in $Y$.
\end{corollary}

\noindent {\bf  Case 2. }  The optimal control
problem {\bf (P)} {\bf with certain control
constraint}.

\medskip

First,  we give  the following hypothesis:

\medskip

\noindent ${\bf (A_{21})}$  Let $f: [0, T]\times
Y\times \widetilde{U}\rightarrow Y$
 and $f^0: [0, T]\times Y\times \widetilde{U}\rightarrow \dbR$ satisfy that
$f$ and $f^0$ are strongly measurable with
respect to $t$ in $(0, T)$, and continuously
Fr\'echet differentiable
 with respect to  $y$ in $Y$ with
$f(t, \cdot, \cdot),  f_y(t, \cdot, \cdot),
f^0(t, \cdot, \cdot)$ and $f^0_y(t, \cdot,
\cdot)$   being continous, respectively.
Moreover,  there exists a  positive constant
$L$,  such  that for any $(t, y, u)\in [0,
T]\times Y\times \widetilde{U}$,
$$
|f_y(t, y, u)|_{\mathcal{L}(Y)}+|f^0_y(t, y,
u)|_{Y'}+|f(t, 0, u)|_{Y}+|f^0(t, 0, u)|\leq L.
$$
\noindent ${\bf (A_{22})}$   For any  $(t, y,
u)\in [0, T]\times  Y\times \widetilde{U}$,
$f(t, y, u)=f_1(t,  y)+B(t)u$  with  $B\in
L^\infty(0, T; \mathcal{L}(U;  Y))$.

\medskip

\noindent  ${\bf  (A_{23})}$  $p=\infty$,
$\widetilde{U}\subseteq U$  is a bounded set and
$\overline{\mbox{co}}\ \widetilde{U}$  has at
least an interior point in $U$.

\medskip

\noindent Under the above  assumptions,  $
\mathcal{U}[0, T]=\Big\{ u\in L^\infty(0, T; U)\
\Big|\ u: (0, T)\rightarrow \widetilde{U}\mbox{
is measurable}  \Big\}, $  and for any $T>0$,
$y(0)\in Y$ and $u(\cdot)\in\mathcal{U}[0, T]$,
$(\ref{61})$ admits a mild solution $y(\cdot)\in
C([0, T]; Y)$ and  $f^0(\cdot, y(\cdot),
u(\cdot))\in L^1(0,  T)$. Also,  the assumptions
$(H_1)$  and $(H_2)$ on Page 130 in \cite{2}
hold.

\medskip

Consider  the following  linear system:
\begin{eqnarray}\label{620}
\left\{
\begin{array}{ll}
\ds\widehat\xi_t(t)=A\widehat\xi(t)+f_{1, y}(t,
\overline{y}(t))\widehat\xi(t)+
B(t)[v(t)-\overline{u}(t)],  &\quad t\in (0, T],\\
\ns\widehat\xi(0)=y^0,
\end{array}
\right.
\end{eqnarray}
where $v\in  \mathcal{U}[0,T]$ and set
$$
M_3=\Big\{  \widehat\xi(T)-y^1\in Y\ \Big|\
\widehat\xi\mbox{ is the solution to }
(\ref{620})\mbox{ with  some  }
v(\cdot)\in\mathcal{U}[0, T]\mbox{  and
}(y^0,y^1)\in S \Big\}.
$$
By Theorem 1.6 on  Page 135 in \cite{2},  if
$M_3$ is {\bf finite  codimensional} in $Y$,
then the optimal pair $(\overline{u}(\cdot),
\overline{y}(\cdot))$ for the  optimal control
problem {\bf (P)} satisfies Pontryagin type
maximum principle, that is, there exists a
nontrivial pair $(\psi^0_2,
\psi_2(\cdot))\in\dbR\times C([0,T]; Y')$,  such
that
\begin{eqnarray*}
\left\{
\begin{array}{ll}
\ds\psi_{2, t}(t)=-A^*\psi_2(t)-\psi^0_2
f^0_y(t,  \overline{y}(t), \overline{u}(t))-
f_{1, y}(t, \overline{y}(t))^*\psi_2(t),
&\quad t\in (0, T],\\
\ns\psi_2(T)\in  Y',
\end{array}\right.
\end{eqnarray*}
and $ H(t, \overline{y}(t), \overline{u}(t),
\psi^0_2,  \psi_2(t))= \max\limits_{u\in
\widetilde{U}} H(t, \overline{y}(t),  u,
\psi^0_2,  \psi_2(t)),\  \mbox{a.e. } t\in
(0,T), $ where $$H(t, y, u, \psi^0, \psi)=\psi^0
f^0(t, y, u)+\langle \psi, f(t, y,
u)\rangle_{Y', Y}, \  \forall\ (t, y, u, \psi^0,
\psi)\in [0, T]\times Y\times
\widetilde{U}\times \dbR\times  Y'.$$

In the rest of this section,  we study  the
finite codimensionality of  $M_3$ with fixed
endpoint constraints. Set $S=\big\{(y^0,
y^1)\big\}$ with $y^0, y^1\in Y$ arbitrarily
given. Then  the solution $\widehat\xi$ to
(\ref{620})  can be rewritten as
$\widehat\xi=z_2+\eta_2$,   where $z_2$ and
$\eta_2$ satisfy, respectively,  that
\begin{eqnarray}\label{621}
\left\{
\begin{array}{ll}
\ds  z_{2,  t}(t)=Az_2(t)+f_{1, y}(t,
\overline{y}(t))z_2(t)+
B(t)v(t),  &\quad t\in (0, T],\\
\ns z_2(0)=0,
\end{array}
\right.
\end{eqnarray}
and
 \begin{eqnarray*}\label{622}
\left\{
\begin{array}{ll}
\ds\eta_{2, t}(t)=A\eta_2(t)+f_{1, y}(t,\overline{y}(t))\eta_2(t)-B(t)\overline{u}(t),  &\quad t\in (0, T],\\
\ns\eta_2(0)=y^0.
\end{array}
\right.
\end{eqnarray*}
Set
\begin{equation}\label{004llz}
M_4=\Big\{ z_2(T)\in Y\ \Big|\ z_2\mbox{ is the
solution to } (\ref{621})\mbox{ with  some  }
v(\cdot)\in\mathcal{U}[0, T] \Big\}.
\end{equation}
Then $M_3=M_4+\big\{\eta_2(T)-y^1\big\}$ and
therefore,   $M_3$  is  finite  codimensional in
$Y$, if and only  if    $M_4$ is  finite
codimensional in $Y$.  Moreover, write
$$\widetilde{\mathcal{U}}=\mathcal{U}[0, T]=\big\{ u\in\mathcal{U}_\infty\ \Big|\ u: (0, T)\rightarrow \widetilde{U}
\mbox{  is measurable}\big\}.$$ By ${\bf
(A_{23})}$,  $\widetilde{\mathcal{U}}$  is a
bounded subset of $\mathcal{U}_\infty$  and
$\overline{\mbox{co}}\ \widetilde{\mathcal{U}}$
has at least an interior point in
$\mathcal{U}_\infty$.

\medskip

By Theorems \ref{t2}-\ref{t3},  one has the
following result on  the finite codimensionality
of $M_4$.

\begin{corollary}  The following  assertions are equivalent:

\medskip

\noindent  {\bf (1)}   The system $(\ref{621})$
is finite  codimensional exactly controllable in
$Y$.

\medskip

\noindent  {\bf (2)} There is a finite
codimensional subspace $\wt Y\subseteq Y'$, such
that  any  solution $\phi$ to the equation
\begin{eqnarray}\label{03llz}
\left\{
\begin{array}{ll}
\ds\phi_t(t)=-A^*\phi(t)-f_{1, y}(t,
\overline{y}(t))^*\phi(t),
&\quad t\in (0, T],\\
\ns\phi(T)=\phi_T,
\end{array}\right.
\end{eqnarray}
satisfies that
\begin{equation}\label{100llz}
|\phi_T|_{Y'}\leq C |B(\cdot)^*\phi|_{L^1(0, T;
U')},\quad\quad\forall  \phi_T\in \wt Y.
\end{equation}

\noindent {\bf (3)}  There  is  a  compact
operator $G$  from  $Y'$ to a Banach space $X$,
such  that  any  solution $\phi$ to
\eqref{03llz} satisfies that
$$
|\phi_T|_{Y'}\leq C\big(|B(\cdot)^*\phi|_{L^1(0,
T;  U')} + |G\phi_T|_{X}\big),\quad\quad\forall
\phi_T\in Y'.
$$

\noindent  {\bf (4)} The  set  $M_4$ $($defined
in $(\ref{004llz}))$ is finite codimensional in
$Y$.
\end{corollary}
\begin{remark}
Similar to  the proofs of $\cite[Theorem
1.1]{zua}$  and Theorem $\ref{t2.1}$ in this
paper, one can  get a weaker criterion  than the
estimate $(\ref{100llz})$. Indeed,
$(\ref{100llz})$ is true,  if and only if  there
is a finite codimensional subspace
$\widetilde{Y}\subseteq Y'$ such that for any
solution  $\phi$ to $(\ref{03llz})$, it holds
that
$$
|\phi_T|_{Y'}\leq C |B(\cdot)^*\phi|_{L^2(0, T;
U')},\quad\quad\forall\ \phi_T\in \wt Y.
$$
\end{remark}

\section{Two examples}\label{sec-exam}

In this section,   two  examples  on
linear  quadratic  control (LQ for short) problems with fixed  endpoint constraints  for wave  and
heat equations   are presented,  respectively.
By the finite codimensional exact controllability and its   equivalent
assertions introduced in this paper,   the finite codimensionality of the sets appeared   in these
optimal control problems
will be verified very easily.

\subsection{Example 1.  An LQ problem for wave equations}\label{sub6.1}

\noindent {\bf 1)  Formulation of  problem}

\medskip

Recall that $\Omega\subset\dbR^n$ is a bounded
domain with a $C^\infty$ boundary $\pa\Omega$.
Put $Q=\Omega\times(0,T)$ and
$\Sigma=\partial\Omega\times(0,T)$. Assume that
$\omega$  is a nonempty open subset  of
$\Omega$. Denote by $\chi_\omega$ the
characteristic function  of $\omega$. Consider
the following wave equation:
\begin{eqnarray}\label{31}\left\{
\begin{array}{ll}
y_{tt}-\Delta y+a(x, t)y=\chi_\omega u  &\mbox{ in } Q,\\
\ns y=0 &\mbox{ on }  \Sigma,\\
\ns y(0)=y_0, \ y_t(0)=y_1  &\mbox{ in }  \Omega,
\end{array}
\right.
\end{eqnarray}
where  $u\in L^2(Q)$ is  the control variable
and $(y, y_t)$ is the state variable,  $(y_0,
y_1)\in H^1_0(\Omega)\times L^2(\Omega)$ is an
initial value, and $a(\cdot) \in L^\infty(Q)$. For a
given $(y^0, y^1)\in H^1_0(\Omega)\times
L^2(\Omega)$,  set
\begin{eqnarray*}
&&\mathcal{U}_{ad}=\Big\{  u(\cdot)\in  L^2(Q)\
\Big|\ \mbox{the   solution }y\mbox{ to
}\eqref{31} \mbox{  satisfies that } (y(T),
y_t(T))=(y^0, y^1) \Big\}
\end{eqnarray*}
and
$$
J(u(\cdot))=\ds\frac{1}{2}\ds\int_Q
\big[q(x, t)|y(x,t)|^2+r(x,t)|u(x,
t)|^2\big]dxdt,
$$
where $y(\cdot)$ is the solution to (\ref{31})
associated to $u(\cdot)$,   and $q,r\in L^\infty(Q)$ are
given functions.

\medskip

Assume that  $(\overline{u}(\cdot),
\overline{y}(\cdot))$  is  an  optimal  pair  of the optimal control problem:
$$J(\overline{u}(\cdot))=\inf\limits_{u(\cdot)\in
\mathcal{U}_{ad} } J(u(\cdot)), $$
 and       set
\begin{eqnarray}\label{80llz}\begin{array}{rl}
&M_5=\Big\{  (y(T),  y_t(T))\in
H^1_0(\Omega)\times L^2(\Omega)\ \Big|\  y\mbox{
is the solution to } (\ref{31})\mbox{ with }
(y_0,  y_1)=(0,  0),\\[3mm]
&\quad\qq\qq\qq\qq\qq\qq\qq\qq\q\mbox{and some }
u(\cdot)\in L^2(Q) \mbox{ satisfying }  |u|_{L^2(Q)}\leq
1 \Big\}.
\end{array}
\end{eqnarray}
Consider the following   equation:
\begin{eqnarray}\label{30llz}\left\{
\begin{array}{ll}
\ds\psi_{tt}-\Delta \psi+a(x,
t)\psi-\psi^0q(x, t)\overline{y}(x, t)=0  &\mbox{
in
}Q,\\
\ns\ds\psi=0 &\mbox{ on }\Sigma,\\
\ns\ds\psi(T)=\psi_1,\  \psi_t(T)=\psi_2
&\mbox{ in }\Omega,
\end{array}
\right.
\end{eqnarray}
where $\psi^0\in\dbR$ and $(\psi_1,  \psi_2)\in
L^2(\Omega)\times  H^{-1}(\Omega)$. It is easy to check that   ${\bf (A_{11})}$ and ${\bf
(A_{12})}$  in Proposition \ref{tt1}  hold.  Hence,    one has  the following
result.
\begin{proposition}\label{90llz}
If  $M_5$  is {\bf finite codimensional} in
$H^1_0(\Omega)\times L^2(\Omega)$,  then there
exist $(\psi^0,  \psi_1,  \psi_2)\in \dbR\times
L^2(\Omega)\times  H^{-1}(\Omega)$, such that
for the corresponding  solution $\psi$ to
$(\ref{30llz})$,  $(\psi^0, \psi(\cdot))\neq (0,
0)$  and
$$\ds
\psi^0r(x, t)\overline{u}(x,
t)+\chi_\omega\psi(x, t)=0,\quad\mbox{a.e. }(x,
t)\in Q.
$$
\end{proposition}

\medskip

\noindent {\bf 2) The geometric control
condition and finite codimensionality}

\medskip

In  this part, under some mild assumptions,  we
prove that  the set  $M_5$  (in Proposition
\ref{90llz})  is finite  codimensional in
$H^1_0(\Omega)\times L^2(\Omega)$,    if and
only if
 the geometric control condition for  $(\O,\o,T)$ in (\ref{31}) holds.
To begin with, let us  recall
some related definitions on  this condition (see \cite{BLR} and also \cite{Dehman-Ervedoza, FYZ, LLTT, LZZ} for more details and related results).

\begin{definition}\label{Def1-0}
${\bf (1)}$  For a wave operator $W=\pa_{tt}
-\Delta $,  a  null bicharacteristic $(t(\cdot),
\hat x(\cdot), \tau(\cdot), p(\cdot)):
\dbR\rightarrow \dbR^{2(n+1)}$ is defined to be
a solution to the  system:
  \bel{rz3}
  \left\{
  \ba{ll}
 t_s(s)=2\tau(s),\  \hat x_{s}(s)=-2p(s),&\\ \ns
 \ds \tau_s(s)=0,\  p_{s}(s)=0,&
  \ea  \right.\ee
with $t(0)=0$, $\tau(0)=1/2$, $\ds \hat
x(0)=\hat x_{0}\in \dbR^n$ and
$p(0)=p_{0}\in\dbR^n$ satisfying that
$|p_{0}|=1/2$.  $(t(\cd),\hat x(\cdot), p(\cd))$
is called a ray of $W$.

\smallskip

%\noindent $(2)$ A ray $\hat x(\cdot)$ of   $W$ is said to start from a domain
%$D\subseteq \dbR^{n}$ at  an initial  time,  if
%$\hat x(0)\in D$.

%\smallskip

%\noindent $(3)$  A ray $\hat x(\cdot)$ of $W$ is said
%to exit $D$ in a finite time, if    there exists
%a
 %time $t_0>0$,  so that the ray $\hat x(\cdot)$ arrives the
%boundary  $\partial D$ of $D$ at  $t_0$ for the   first
%time, i.e., $\hat x(t_0)\in \partial D$.

%\smallskip

\noindent ${\bf  (2)}$  For any $T>0$  and open
set $\Omega$  of $\dbR^n$,   $(t(\cd),\hat
x(\cdot), p(\cdot)):  [0, T]\rightarrow
\overline{\Omega}\times\dbR^{n}$  with $\hat
x(0), \hat x(T)\in \Omega$ is called a
generalized ray of  $W$ in $\overline{\Omega}$,
if there exists a  partition
$0=s_0<s_1<\cdots<s_m=T$ for an  $m\in\dbN$,
such that   for any $j=0, 1, \cdots, m-1$,
$(\hat
 x(s), p(s))\big|_{s_j\leq s\leq s_
{j+1}}=(\hat x^j(s), p^j(s))$ satisfies
$(\ref{rz3})$, $\hat x(s_{k})\in\partial \Omega$
$(k=1, 2, \cdots, m-1)$, and  the following law
of geometric optics holds:
%\begin{equation}\label{op5}
%\begin{array}{ll}\ds
$$\ds p^{k+1}(s_{k})=p^k(s_{k})-2 \nu(\hat
x^k(s_{k}))^\top  p^k(s_{k}) \nu(\hat
x^k(s_{k})),
$$
%\end{array}
%\end{equation}
where   $\nu(x)$ denotes  the unit outer normal vector
on  $x\in \partial \Omega$.  $s_{k}$   is
called  the $k$-th reflected instant of this
generalized ray.  A  generalized
ray is denoted by $\big\{(t,\hat x^j(t),
p^j(t))\;\big|\; t\in[s_j,
s_{j+1}]\big\}_{j=0}^{m-1}$.

\smallskip

\noindent ${\bf (3)}$ $(\O,\o,  T)$ in
$(\ref{31})$ is  called to satisfy   the
geometric control  condition, if for any
generalized ray $\big\{(t(s),\hat x^j(s),
p^j(s))\;\big|\; s\in[s_j,
s_{j+1}]\big\}_{j=0}^{m-1}$ of  $W$ in
$\overline{\O}$, there are a $j\in\{0, 1,
\cdots, m-1\}$  and $s_0\in[s_j, s_{j+1}]$, such
that $\hat x^j(s_0)\in\omega$.

\end{definition}

\begin{remark}
It is easy to show that a null bicharacteristics
of $W$  is a straight line in $\dbR^{2(n+1)}$.
Noting that $t(s)\equiv s$, in the rest of this
paper, we simply denote by $\big\{(t,\hat x(t),
p(t))\;\big|\; t\in [0,T]\big\}$ the generalized
ray.
\end{remark}
%
%\begin{remark}
%
%By
%$(\ref{rz3})$,   $(\ref{op5})$  may be
% reduced to  the condition
%
%\begin{equation}\label{ggop5}
%x_t^{k+1}(s_{k})= x_t^k(s_{k})-2
%\nu(x^k(s_{k}))^\top x_t^k(s_{k})
%\nu(x^k(s_{k})).
%\end{equation}
%
%This means that the direction of $x^{k+1}(\cdot)$ at
%$t=s_{k}$ is obtained from that of  $x^k(s_k)$  by a reflection with
%respect to $\nu(x^k(s_{k}))$. \end{remark}

\medskip

In order to prove that $M_5$ is   finite  codimensional in $H^1_0(\Omega)\times L^2(\Omega)$ (in Proposition \ref{90llz}),
consider the following backward wave equation:
\begin{eqnarray}\label{32}\left\{
\begin{array}{lll}
\ds\phi_{tt}-\Delta \phi+a(x,
t)\phi=0  &\mbox{ in
}Q,\\
\ns\ds\phi=0 &\mbox{ on }\Sigma,\\
\ns\ds\phi(T)=\phi_1,\  \phi_t(T)=\phi_2
&\mbox{ in }\Omega,
\end{array}
\right.
\end{eqnarray}
where $(\phi_1, \phi_2)\in L^2(\Omega)\times
H^{-1}(\Omega)$.  By Corollary \ref{LLLL}, we
have the following result.
\begin{proposition}\label{llz11}
If  $(\Omega, \omega,  T)$   satisfies the
geometric control condition,   then $M_5$
$($defined in $(\ref{80llz}))$ is finite
codimensional  in $H^1_0(\Omega)\times
L^2(\Omega)$.
\end{proposition}

\noindent {\bf Proof. } First,  denote by $\phi$
the solution to (\ref{32}) and  set
$\phi=\xi_1+\xi_2$, where $\xi_j$ $(j=1, 2)$
satisfy
 \begin{eqnarray*}\left\{
\begin{array}{lll}
\ds\xi_{1, tt}-\Delta \xi_1=0  &\mbox{ in
}Q,\\
\ns\ds\xi_1=0 &\mbox{ on }\Sigma,\\
\ns\ds\xi_1(T)=\phi_1,\  \xi_{1, t}(T)=\phi_2
&\mbox{ in }\Omega,
\end{array}
\right.
\end{eqnarray*}
and
\begin{eqnarray}\label{8.17-eq1}
\left\{
\begin{array}{lll}
\ds\xi_{2, tt}-\Delta \xi_2+a(x,
t)\phi=0  &\mbox{ in
}Q,\\
\ns\ds\xi_2=0 &\mbox{ on }\Sigma,\\
\ns\ds\xi_2(T)=0,\  \xi_{2, t}(T)=0  &\mbox{ in
}\Omega.
\end{array}
\right.
\end{eqnarray}
Define  operators $G_1$, $G_2$ and $G_3$ as
follows:
$$
G_1: L^2(\Omega)\times H^{-1}(\Omega) \to
L^2(Q),\q \ds G_1(\phi_1, \phi_2)=\phi,\q
\forall (\phi_1, \phi_2)\in L^2(\Omega)\times
H^{-1}(\Omega),
$$
where $\phi$ is the solution to \eqref{32} with
the initial value $(\phi_1, \phi_2)$;
$$
G_2: L^2(Q)  \to L^\infty(0,T;H_0^1(\Omega))\cap
W^{1,\infty}(0,T;L^2(\Omega)),\q
\ds G_2(\phi)=\xi_2,
$$
where $\xi_2$ is the solution to
\eqref{8.17-eq1} associated to  $\phi\in
L^2(Q)$; and
$$G_3: L^\infty(0,T;H_0^1(\Omega))\cap
W^{1,\infty}(0,T;L^2(\Omega))\to L^2(Q),\q
\ds G_3(\xi_2)=\xi_2.
$$
Let $G=G_3G_2G_1$. Then
$
G: L^2(\Omega)\times H^{-1}(\Omega)\to L^2(Q),\
\ds G(\phi_1, \phi_2)=\xi_2.
$

\smallskip

By the well-posedness results of wave equations,
$G$ is a compact operator.  Since $(\Omega,
\omega,  T)$  fulfills the geometric control
condition,  by \cite{BLR1}, we have that
$$
\ds|(\phi_1, \phi_2)|_{L^2(\Omega)\times
H^{-1}(\Omega)}^2 \ds\leq C\int^T_0\int_\omega
\xi_{1}^2dxdt,\q\forall\  (\phi_1, \phi_2)\in L^2(\Omega)\times
H^{-1}(\Omega).
$$
It follows that
$$
\begin{array}{ll}
&\ds|(\phi_1, \phi_2)|_{L^2(\Omega)\times
H^{-1}(\Omega)}^2 \ds \leq C\int^T_0\int_\omega
(\phi^2+\xi^2_2)dxdt\leq C|\chi_{\omega} \phi|^2_{L^2(Q)}+C|G(\phi_1,
\phi_2)|_{L^2(Q)}^2.
\end{array}
$$
Take  $Y'=L^2(\Omega)\times H^{-1}(\Omega)$ and
$X=L^2(Q)$ in  (\ref{20llz}), and   note that
$f_u(\cdot, \overline{y}(\cdot),
\overline{u}(\cdot))^*\phi=\chi_{\omega} \phi$,
by Corollary \ref{LLLL}, $M_5$ $($defined in
$(\ref{80llz}))$ is finite codimensional  in
$H^1_0(\Omega)\times L^2(\Omega)$.
 \endpf

\begin{remark}
By Theorem $\ref{t2}$,    the result in
Proposition $\ref{llz11}$ implies  that   under
the geometric control condition,  the wave
equation $(\ref{31})$  is finite codimensional
exactly controllable. Note that under the same
condition,   the exact controllability of this
wave equation is still open. However,  it is
rather easy to show the  weaker finite
codimensional controllability, while this
controllability  is enough for us  to study the
optimal control problem with constraints.
\end{remark}

Next,  we prove that, under some mild
assumptions,  if the geometric control condition
fails, the system $(\ref{31})$ is not finite
codimensional exactly controllable any more.
Hence, the finite codimensionality of $M_5$
fails in this  case.

\begin{proposition}\label{8.18-th1}
Assume that  there is a generalized ray
$\big\{(t, \hat{x}^j(t), p^j(t))\;\big|\;
t\in[s_j,s_{j+1}]\big\}_{j=0}^{m-1}$ of   $W$
in $\cl{\O}$,  such that
\medskip

\noindent {\bf (1)} it does not meet $\o$, i.e.,
$x^{j}(t)\not\in\cl{\o}$,  for any
$j\in\{0,1,\cdots,m-1\}$ and $t\in
[s_j, s_{j+1}]$; and

\medskip

\noindent {\bf (2)}   it always meets $\partial\Omega$ transversally, i.e.,
  $
\nu(\hat{x}^k(s_{k}))^\top  p^k(s_{k})\neq  0,\  \forall\; k\in\{1, 2,\cdots, m-1\},
 $
where $ s_{k}$ is the $k$-th reflected
instant of this generalized ray.

\medskip

\noindent Then the system \eqref{31} is not
finite  codimensional exactly controllable.
\end{proposition}

By Theorem \ref{t2.1},  the finite codimensional
exact controllability of  the system (\ref{31})
may be reduced to a suitable observability
estimate  (\ref{221}) for solutions $\phi$ to
(\ref{32}). Let $\hat\phi(x, t)=\phi(x, T-t)$.
Then it is easy to check that $\hat\phi(\cdot)$
solves
\begin{eqnarray}\label{5.8-eq5}
\left\{
\begin{array}{lll}
\ds\ds\hat\phi_{tt}-\Delta \hat\phi+\hat
a(x, t)\hat\phi=0 &\mbox{ in }Q,\\\ns\ds
\hat\phi=0 &\mbox{ on }\Sigma,\\\ns\ds \hat
\phi(0)=\phi_1,\  \hat \phi_t(0)=\phi_2
&\mbox{ in }\Omega,
\end{array}
\right.
\end{eqnarray}
where $\hat a(x, t)=a(x, T-t)$ and $(\phi_1, \phi_2)\in L^2(\Omega)\times H^{-1}(\Omega)$.
Also,  the   observability estimate for  solutions $\phi$ to  (\ref{32})  is equivalent to the following one for   solutions $\hat\phi$  to
(\ref{5.8-eq5}):
\begin{equation}\label{llz101}
|(\phi_1,  \phi_2)|_{L^2(\Omega)\times H^{-1}(\Omega)}\leq
C|\hat\phi|_{L^2(0, T; L^2(\omega))}, \quad\forall \  (\phi_1,  \phi_2)\in \widehat{M},
\end{equation}
where $\widehat{M}$ is a finite  codimensional subspace of $L^2(\Omega)\times H^{-1}(\Omega)$.

\medskip

As preliminaries  to prove Proposition
\ref{8.18-th1}, some lemmas  are given  in
order. First, consider the following wave
equation:
\begin{eqnarray}\label{5.8-eq211}
\left\{
\begin{array}{lll}
\ds\f_{tt}-\Delta \f + \int^t_0 \hat a(x,
s)\f_s(x, s)ds=0  &\mbox{ in }Q,\\
\ns\ds\f=0 &\mbox{ on }\Sigma,\\
\ns\ds\f(0)=\f_1,\ \f_t(0)=\f_2  &\mbox{ in
}\Omega,
\end{array}
\right.
\end{eqnarray}
where $(\f_1,\f_2)\in H_0^1(\O)\times L^2(\O)$.
Then the observability estimate   (\ref{llz101})  for solutions $\hat\phi$ to  (\ref{5.8-eq5}) implies  a suitable  estimate
for solutions $\f$ to (\ref{5.8-eq211}) (Notice that conversely, it  may  be untrue).
\begin{lemma}\label{5.8-lm1}
If any solution  $\hat\phi$ to the  equation
\eqref{5.8-eq5} satisfies $(\ref{llz101})$,
then there exists a finite codimensional
subspace $M_6$ of $H_0^1(\O)\times L^2(\O)$,
such that  solutions $\f$ to
$(\ref{5.8-eq211})$ satisfy
\begin{equation}\label{5.8-lm1-eq1}
|(\f_1,\f_2)|_{H_0^1(\O)\times L^2(\O)}\leq
C|\f_t|_{L^2(0,T;L^2(\o))},  \q\forall\
(\f_1,\f_2)\in M_6.
\end{equation}
\end{lemma}
\noindent {\bf Proof.}  For arbitrarily given
$(\f_1,\f_2)\in H_0^1(\O)\times L^2(\O)$, let
$\f$ be the corresponding   solution to
(\ref{5.8-eq211}). Set $\hat\phi=\f_t$. Then
$\hat\phi$ solves \eqref{5.8-eq5} with the
initial value $(\hat\phi(0),
\hat\phi_t(0))=(\f_2,\D\f_1)\in
L^2(\Omega)\times H^{-1}(\Omega)$.  By the
assumption,
 there exists a finite
dimensional subspace $M_6$ of $L^2(\O)\times
H^{-1}(\O)$,  such that $L^2(\O)\times
H^{-1}(\O)=M_6\oplus\widehat{M}$, and
 for any $(\f_2,\D\f_1)\in
\widehat{M}$, it holds that
\begin{equation}\label{5.8-eq3}
|(\f_2,\D\f_1)|_{L^2(\O)\times H^{-1}(\O)}\leq
C|\f_t|_{L^2(0,T;L^2(\o))}.
\end{equation}
Denote by  $\mathcal{A}$ be the Laplacian
operator with homogeneous  Dirichlet boundary
condition.  Set
$$
M_7=\Big\{(\f_1,\f_2)\in H_0^1(\O)\times
L^2(\O)\ \|\   \f_2=\hat \phi_1,   \f_1 =
(-\mathcal{A})^{-1}\hat\phi_2 \mbox{ for some }
(\hat \phi_1,\hat\phi_2)\in \widehat{M} \Big\}.
$$
and
$$
M_8=\Big\{(\f_1,\f_2)\in H_0^1(\O)\times
L^2(\O)\  \|\ \f_2=\hat \phi_1, \; \f_1 =
(-\mathcal{A})^{-1}\hat\phi_2 \mbox{ for some }
(\hat \phi_1,\hat\phi_2)\in  M_6 \Big\}
$$
Then by (\ref{5.8-eq3}),  it is clear that
(\ref{5.8-lm1-eq1})  holds. Also,  $M_8$ is
finite dimensional  and $H^1_0(\O)\times
L^2(\O)=M_7\oplus M_8$.  The proof is completed.
\endpf

\medskip

\medskip

Further,  a priori estimate  for  a compact operator  is  presented.
\begin{lemma}\label{5.8-lm2}
Assume that $G$ is a compact operator from
$H_0^{1}(\O)\times L^2(\O)$ to itself. Then for
any $\d>0$, there is a positive  constant $C$
such that
\begin{equation}\label{5.8-lm2-eq1}
\begin{array}{ll}\ds
&\displaystyle|G(\f_1,\f_2)|_{H_0^{1}(\O)\times L^2(\O)}\\[4mm]
&\displaystyle\leq
\d|(\f_1,\f_2)|_{H_0^{1}(\O)\times L^2(\O)}+
C|(\f_1,\f_2)|_{L^2(\O)\times H^{-1}(\O)},\quad \forall \   (\f_1,\f_2)\in
H_0^{1}(\O)\times L^2(\O).
\end{array}
\end{equation}
\end{lemma}
\noindent{\bf Proof.} Assume that
(\ref{5.8-lm2-eq1}) fails. Then there exist
$\delta_0>0$  and a sequence
$\{(\f_{1,n},\f_{2,n})\}_{n=1}^\infty$  in $
H_0^{1}(\O)\times L^2(\O)$, such that for any
$n\in\dbN$,
\begin{eqnarray}\label{8.17-eq2}
\left\{
\begin{array}{ll} \ds
\big|G(\f_{1,n},\f_{2,n})\big|_{H_0^{1}(\O)\times
L^2(\O)}=1,&\\
\ns\ds
\big|(\f_{1,n},\f_{2,n})\big|_{H_0^{1}(\O)\times
L^2(\O)}\leq\frac{1}{\d_0}, &\\
\ns\ds |(\f_{1,n},\f_{2,n})|_{L^2(\O)\times
H^{-1}(\O)}\leq \frac{1}{n}.&
\end{array}\right.
\end{eqnarray}
Hence,  there exist a subsequence of  $
\{(\f_{1,n},\f_{2,n})\}_{n=1}^\infty$ (still
denoted  by  itself) and $(\f_{1},\f_{2})\in
H_0^{1}(\O)\times L^2(\O)$, such that
$$
(\f_{1, n},\f_{2, n})\rightarrow (\f_{1},\f_{2})
\mbox{  weakly in  }H_0^{1}(\O)\times L^2(\O)
\mbox{ as } n\rightarrow +\infty.
$$
Since the embedding from $H_0^{1}(\O)\times
L^2(\O)$ to $L^2(\O)\times H^{-1}(\O)$ is
compact, then
$$ \lim_{n\to\infty}(\f_{1, n},\f_{2, n})=
(\f_{1},\f_{2}) \mbox{  weakly in }L^2(\O)\times
H^{-1}(\O). $$
This, together with the third inequality in
\eqref{8.17-eq2}, deduces that
$(\f_{1},\f_{2})=(0,0)$. Noting that $G$ is
compact, we obtain that
$$
\lim\limits_{n\to\infty} G(\f_{1,  n},\f_{2, n})
=G(\f_{1},\f_{2}) \mbox{ in }H_0^{1}(\O)\times
L^2(\O),
$$
which implies that   $
|G(\f_{1},\f_{2})|_{H_0^{1}(\O)\times
L^2(\O)}=1$. It contradicts  that   $(\f_1,
\f_2)=(0, 0)$. Hence,  (\ref{5.8-lm2-eq1})
holds.
\endpf

\medskip

Further, we construct a family of  solutions to
the equation  (\ref{5.8-eq211}).
\begin{lemma}\label{8.18-lm1}
Suppose that all assumptions in Proposition
$\ref{8.18-th1}$ hold  and $\hat a\in
L^\infty(0,  T;  W^{1, \infty}(\Omega))$. Then
there exist a family of solutions
$\{\f_\e\}_{\e>0}$ to \eqref{5.8-eq211}, and
positive   constants $c_1$ and $c_2$,
independent of $|\nabla\hat a|_{L^\infty(Q)}$,
such that for    any $0<\e<1$,
\begin{equation}\label{or17}\left\{
\begin{array}{ll}\ds
|\f_{\e,t}(\cd,  0)|_{L^2(\O)}\ge c_1,&\\
\ns\ds |\f_{\e, t}(\cdot, 0)|_{L^2(\O)}\leq
c_2,\ |\f_\e(\cdot, 0)|_{H^1_0(\O)}\leq
c_2,&\\\ns\ds
|\f_\e|_{H^1(0,T;\,L^2(\o))}=(|\hat a|_{L^\infty(0,  T;  W^{1, \infty}(\Omega))}+1)\mathcal O(\e^{1/2}),&\\
\ns\ds
|\f_{\e}(\cd,  0)|_{L^2(\O)}+|\f_{\e,t}(\cd,  0)|_{H^{-1}(\O)}=\mathcal O(\e^{1/2}),&
\end{array}\right.
\end{equation}
where $\mathcal O(\e^{\alpha})$  denotes a function of order $\e^{\alpha}$ for $\alpha>0$.
\end{lemma}
The proof of  Lemma \ref{8.18-lm1} is similar to
that of Theorem 7.1 in \cite{FYZ}.  We shall only give its
sketch  in Appendix for completeness.

\medskip

Now, we are in a position  to prove Proposition \ref{8.18-th1}.

\medskip

\noindent {\bf Proof of Proposition
\ref{8.18-th1}.} Assume that the system
\eqref{31} is finite codimensional exactly
controllable.  Then by  Theorem \ref{t2.1} and
Lemma \ref{5.8-lm1}, there are a finite
dimensional subspace $M_8$ and finite
codimensional subspace $M_7$ of  $
H^{1}_0(\O)\times L^2(\O)$, such that
$H^1_0(\O)\times L^2(\O)=M_8\oplus M_7$,  and
for any $(\f_1, \f_2)\in  M_7$, the
corresponding solution $\f$ to the   equation
\eqref{5.8-eq211} satisfies
$(\ref{5.8-lm1-eq1})$. Denote by $\mathbb
P_{M_8}$ the projection   from $H^1_0(\O)\times
L^2(\O)$ to $ M_8$.  Then for any $(\f_1,
\f_2)\in H^1_0(\O)\times L^2(\O)$,  it holds
that
\begin{eqnarray}\label{5.8-eq9}
\begin{array}{ll}\ds
|(\f_1,\f_2)|^2_{H^1_0(\O)\times
L^2(\O)}\\
\ns\ds\ds=|(\f_1,  \f_2)- \mathbb P_{
M_8}(\f_1,\f_2)|^2_{H^1_0(\O)\times
L^2(\O)}+|\mathbb P_{
M_8}(\f_1,\f_2)|^2_{H^1_0(\O)\times L^2(\O)}\\
\ns \ds\leq
C|\hat\f_t|^2_{L^2(0,T;L^2(\o))}+|\mathbb P_{
M_8}(\f_1,\f_2)|^2_{H^1_0(\O)\times L^2(\O)}\\
\ns \ds \leq C|\f_t|^2_{L^2(0,T;L^2(\o))}+
C|\check\f_t|^2_{L^2(0,T;L^2(\o))} +|\mathbb
P_{M_8}(\f_1,\f_2)|^2_{H^1_0(\O)\times
L^2(\O)}\\ \ns \ds\leq C\big(|
\f_t|^2_{L^2(0,T;L^2(\o))} +|\mathbb P_{
M_8}(\f_1,\f_2)|^2_{H^1_0(\O)\times
L^2(\O)}\big),
\end{array}
\end{eqnarray}
where $\hat\f$ and $\check\f$ are the solutions
to the equation  \eqref{5.8-eq211},
respectively, with the initial values $(\f_1,
\f_2)- \mathbb P_{ M_8}(\f_1,\f_2)$  and
$\mathbb P_{M_8}(\f_1,\f_2)$.

\smallskip

Since $\mathbb  P_{M_8}$ is compact,  by Lemma
\ref{5.8-lm2}, for any $\rho>0$,
$$
|\mathbb P_{M_8}(\f_1,\f_2)|^2_{H^1_0(\O)\times
L^2(\O)} \leq
\rho|(\f_1,\f_2)|^2_{H^1_0(\O)\times
L^2(\O)} +
C|(\f_1,\f_2)|^2_{L^2(\O)\times
H^{-1}(\O)}.
$$
This, together with \eqref{5.8-eq9}, implies
that
\begin{equation}\label{5.8-eq10}
\begin{array}{ll}\ds
&|(\f_1,\f_2)|^2_{H^1_0(\O)\times L^2(\O)}
\ds\leq C|\f_{t}|^2_{L^2(0,T;L^2(\o))}
+C|(\f_1,\f_2)|^2_{L^2(\O)\times
H^{-1}(\O)}.
\end{array}
\end{equation}

For  $\hat a(\cdot)\in L^\infty(Q)$  and  any
$\rho>0$,  there exists a  function
$a_\rho(\cdot)\in L^\infty(0,  T;  W^{1,
\infty}(\Omega))$,  such that $ |a_\rho-\hat
a|_{L^2(Q)}<\rho$ and
$|a_\rho|_{L^\infty(Q)}\leq C$, where $C$ is
independent of $\rho$.   For any
$\varepsilon>0$,  by Lemma \ref{8.18-lm1}, there
exist $(\f^\rho_{\e,  1}, \f^\rho_{\e, 2})\in
H^1_0(\O)\times L^2(\O)$,  such that  the
solution $\f_\e^\rho(\cdot)$ to
(\ref{5.8-eq211}) with $\hat
a(\cdot)=a_\rho(\cdot)$ and $(\f_{1},
\f_2)=(\f^\rho_{\e, 1}, \f^\rho_{\e, 2})$
satisfies the  conclusions  in (\ref{or17}).
More precisely,
$$\left\{
\begin{array}{ll}\ds
|\f^\rho_{\e,t}(\cd,  0)|_{L^2(\O)}\ge c_1,&\\
\ns\ds |\f^\rho_{\e, t}(\cdot, 0)|_{L^2(\O)}\leq
c_2,\ |\f^\rho_\e(\cdot, 0)|_{H^1_0(\O)}\leq
c_2,&\\\ns\ds |\f^\rho_\e|_{H^1(0,T;\,L^2(\o))}=
(|a_\rho|_{L^\infty(0,  T;  W^{1, \infty}(\Omega))}+1)\mathcal O(\e^{1/2}),&\\
\ns\ds
|\f^\rho_{\e}(\cd,  0)|_{L^2(\O)}+|\f^\rho_{\e,t}(\cd,  0)|_{H^{-1}(\O)}=\mathcal O(\e^{1/2}),&
\end{array}\right.
$$
where $c_1$ and $c_2$ are  independent of $\rho$
and $\e$. Denote by $\psi^\rho_\e(\cdot)$  the
solution  to  (\ref{5.8-eq211}) with $\hat
a(\cdot)\in L^\infty(Q)$ and   $(\f_{1},
\f_2)=(\f^\rho_{\e, 1}, \f^\rho_{\e, 2})$. Set
$z_\e^\rho(\cdot)=\f^\rho_\e(\cdot)-\psi^\rho_\e(\cdot)$.
Then $z_\e^\rho(\cdot)$ satisfies
\begin{eqnarray*}
\left\{
\begin{array}{lll}
\ds z^\rho_{\e, tt}-\Delta z^\rho_\e + \int^t_0 a_\rho(x,
s)  z^\rho_{\e,  s}(x, s)ds= \int^t_0  [\tilde a(x,
s)-a_{\rho}(x, s)]  \psi^\rho_{\e,  s}(x, s)ds  &\mbox{ in }Q,\\
\ns\ds  z^\rho_\e=0 &\mbox{ on }\Sigma,\\
\ns\ds z^\rho_\e(0)=0,\  z^\rho_{\e, t}(0)=0 &\mbox{ in
}\Omega,
\end{array}
\right.
\end{eqnarray*}
and
%\begin{eqnarray*}
%\begin{array}{ll}
$$\ds |z^\rho_\e|_{H^1(0, T; L^2(\o))}\leq  C \rho |\psi^\rho_{\e, t}|_{L^2(Q)}
%&
%\\\ns\ds
\leq C\rho (|\f^\rho_{\e,
1}|_{H^1_0(\O)}+|\f^\rho_{\e, 2}|_{L^2(\O)})
\leq C c_2 \rho.
%&
%\end{array}
%\end{eqnarray*}
$$Hence,
$$
|\psi^\rho_\e|_{H^1(0, T; L^2(\o))}\leq C c_2
\rho+C(|a_\rho|_{L^\infty(0,  T;  W^{1,
\infty}(\Omega))}+1)\e^{1/2}.
$$
By   (\ref{5.8-eq10}), it follows that
$$
c_1\leq C c_2 \rho+C(|a_\rho|_{L^\infty(0,  T;
W^{1, \infty}(\Omega))}+1)\e^{1/2}+C\e^{1/2}.
$$
Take $\rho>0$ sufficiently small, such that $C
c_2 \rho<c_1/2$. Then, choose  $\e>0$ small
enough, such that $C(|a_\rho|_{L^\infty(0,  T;
W^{1,
\infty}(\Omega))}+1)\e^{1/2}+C\e^{1/2}<c_1/2$.
This leads to a contradiction.  Hence,   the
system \eqref{31} is not finite codimensional
exactly controllable.
\endpf

\subsection{Example 2. An LQ problem for heat  equations}

Consider the following heat equation:
\begin{eqnarray}\label{18llz}\left\{
\begin{array}{lll}
\ds y_{t}-\Delta y=\chi_\omega u  &\mbox{ in
}Q,\\
\ns\ds y=0 &\mbox{ on }\Sigma,\\
\ns\ds y(0)=y_0  &\mbox{ in }\Omega,
\end{array}
\right.
\end{eqnarray}
where  $u\in L^2(Q)$ is  the control variable,
$y$ is the state variable and $y_0\in
L^2(\Omega)$ is an initial value. Set
\begin{eqnarray*}
&&\mathcal{U}_{ad}=\big\{  u(\cdot)\in  L^2(Q)\
\big|\ \mbox{the   solution }y\mbox{ of
}(\ref{18llz}) \mbox{  satisfies that } y(T)=0
\big\}
\end{eqnarray*}
and
$$
J(u(\cdot))=\ds\frac{1}{2}\ds\int_Q
\big[q(x, t)|y(x,t)|^2+r(x,t)|u(x,
t)|^2\big]dxdt,
$$
where $y(\cdot)$ is the solution  to  (\ref{18llz})
associated to $u(\cdot)$,   and $q,r\in L^\infty(Q)$ are
given  functions.
Assume that  $(\overline{u}(\cdot),
\overline{y}(\cdot))$  is  an  optimal  pair  of the optimal control problem:
$$J(\overline{u}(\cdot))=\inf\limits_{u(\cdot)\in
\mathcal{U}_{ad} } J(u(\cdot)).$$
Write
$$
M_9=\Big\{  y(T)\in L^2(\Omega)\ \Big|\ y\mbox{
is the solution to } (\ref{18llz})\mbox{  with }
y_0=0 \mbox{ and  }v\in L^2(Q): |v|_{L^2(Q)}\leq
1 \Big\}.
$$
Then  similar to  the argument of Proposition \ref{90llz}, in order to guarantee (\ref{0llz}) in
Pontryagin type  maximum  principle to hold, it is required to check if  the set $M_9$ is
finite codimensional.
By Corollary \ref{LLLL} and a contradiction argument,   it is  easy to check  the following negative result  on the finite codimensionality  for the heat equation.
\begin{proposition}\label{pllz}
For any $\Omega$, $\omega$ and $T>0$, $M_9$ is not finite
codimensional in $
L^2(\Omega)$.
\end{proposition}
\noindent {\bf Proof. } By Corollary \ref{LLLL},  it suffices  to prove  that $(\ref{20llz})$  fails  for  the following heat equation:
\begin{eqnarray}\label{21llz}\left\{\!\!\!
\begin{array}{lll}
&\phi_{t}+\Delta \phi=0  &\mbox{ in
}Q,\\
\ns&\phi=0 &\mbox{ on }\Sigma,\\
\ns&\phi(T)=\phi_T  &\mbox{ in }\Omega,
\end{array}
\right.
\end{eqnarray}
with $\phi_T\in  L^2(\Omega)$. In
(\ref{20llz}), take  $U'=Y'=L^2(\Omega)$  and
notice that $f_u(\cdot, \overline{y}(\cdot),
\overline{u}(\cdot))^*\phi=\chi_{\omega} \phi$.

\smallskip

Assume  that $\{\phi^j\}_{j=1}^\infty$  are the  solutions  to    $(\ref{21llz})$
corresponding to  terminal values $\{\phi_T^j\}_{j=1}^{\infty}\subseteq L^2(\Omega)$ with $|\phi_T^j|_{L^2(\Omega)}=1$. Then  there  are a  subsequence of $\{\phi_T^j\}_{j=1}^{\infty}$
$($still  denoted by itself$)$ and $\phi_T\in L^2(\Omega)$,
such that   as  $j\rightarrow  \infty$,
\begin{equation}\label{*LL*}
\phi_{T}^j  \rightarrow  \phi_T\  \mbox{
weakly}\mbox{ in } L^2(\Omega),\q
\phi^j\rightarrow  \phi \  \mbox{ in }
L^2(Q)\q\mbox{ and }\q G(\phi^j_T)\rightarrow
G(\phi_T) \ \mbox{ in }X,\end{equation} where
$\phi$ is  the solution to $(\ref{21llz})$
associated  to $\phi_T$  and $G$  is a compact
operator from  $L^2(\Omega)$  to a Banach  space
$X$. If  $(\ref{20llz})$  is true,  then for any
$j\in\dbN$,
$$|\phi_T^j|_{L^2(\Omega)}\leq C\big(|\phi^j|_{L^2(Q)} +
|G\phi_T^j|_{X}\big).$$ By (\ref{*LL*}), this implies
$\lim\limits_{j\rightarrow
\infty}\phi_{T}^j=\phi_T$ in $L^2(\Omega)$, which leads to a contradiction. Therefore, the  set $M_9$
is  not  finite codimensional in $L^2(\Omega)$.
\endpf

\medskip

Proposition   \ref{pllz}  indicates that  the
finite codimensionality of $M_9$ fails for an LQ
problem of heat equations with fixed endpoint
constraints. Hence, the non-triviality of the
Lagrange type multiplier cannot be established by the
method introduced in this paper.

%%%%%%%%%%%%%%%%%%%%%%%%%%%%%%%%%%%%%%%%%%%%%%%%%%%%%%%%%%%

%%%%%%%%%%%%%%%%%%%%%%%%%%%%%%%%%%%%%%%%%%%%%%%%%%%%%%%%%%%

\section{Appendix}

This section is devoted to  proving  Lemma \ref{8.18-lm1}.

\medskip

\noindent {\bf Proof  of Lemma \ref{8.18-lm1}.} We shall borrow some ideas from
\cite{FYZ, M, R}.
The proof is divided into four steps.

\smallskip

\noindent {\bf Step 1.}  We construct highly
concentrated approximate solutions to
the  hyperbolic equation:
\begin{equation}\label{8.18-eq1}
\f_{tt}-\Delta \f+
F(\f) =0 \qq\mbox{ in }\dbR^n\times
[-T, T],
\end{equation}
where $\ds F(\f) = \int^t_0 \hat a(x, s)\f_s(x,
s)ds $ and   $\hat a\in L^\infty(0, T;  W^{1, \infty}(\dbR^n))$.

Given a generalized ray $(\hat x(\cdot), p(\cdot))$ of
the wave operator $W$,  for any   $\e\in (0,1)$,
construct a family of
approximate solutions  $\phi_\e$ to  the equation
\eqref{8.18-eq1} as
\begin{equation}\label{ssrz5}
\phi_\e(x, t)=\e^{1-\frac{n}{4}}
c(t)e^{{i\psi(x, t)}/\e}+
\e^{2-\frac{n}{4}}\int^t_0 A(s) e^{{i\psi(x,
s)}/\e}ds,
\end{equation}
with
%
%\begin{equation}\label{rz6}
$$\ds
\psi(x,  t)=p^\top(t)[x-\hat
x(t)]+\frac{1}{2}\big[x-\hat x(t)\big]^\top
M(t)\big[x-\hat x(t)\big],
$$
%\end{equation}
%
where  $M(t)$ is a complex symmetric
matrix  with positive
definite imaginary part. The construction of
approximate solutions $\phi_\e$ means  an
appropriate choice of $c(\cdot)$, $M(\cdot)$ and $A(\cdot)$.

Notice that in the subsequent estimates,  we are only concerned with
the dependence of constants on $|\nabla \hat a|_{L^\infty(Q)}$,
rather than $|\hat a|_{L^\infty(Q)}$.

By (\ref{ssrz5}) and a direct computation, it is
easy to check that
\begin{equation}\label{Wu}
\phi_{\e, tt}-\Delta \phi_\e+ F(\phi_\e)
=\e^{2-\frac{n}{4}}r_1+ \e^{1-\frac{n}{4}}r_2
+\e^{-\frac{n}{4}}r_3 + \e^{-1-\frac{n}{4}}r_4,
\end{equation}
where
\begin{eqnarray}\label{r1}\left\{
\begin{array}{ll}
\ds r_1=A_t(t)e^{{i\psi(x,   t)}/\e}+\int^t_0 \hat{a}(x, s)A(s)e^{{i\psi(x,   s)}/\e} ds,&\\\ns\ds
\ds r_2\ds=c_{tt}(t) e^{{i\psi(x, t)}/\e}+i\psi_t(x, t)A(t)e^{{i\psi(x,   t)}/\e}&\\\ns\ds
\q\q
+ \int_0^t \hat a(x,  s)c_s(s)
e^{{i\psi(x,  s)}/\e}ds-
i\int^t_0 A(s)\Delta \psi(x, s)e^{{i\psi(x,   s)}/\e}ds,&\\  \ns\ds
r_3=i\big[2c_t(t)\psi_t(x,  t) + c(t)(W\psi)(x,  t)
\big]e^{{i\psi(x,   t)}/\e}&\\\ns\ds
\quad\q+\int^t_0 A(s)|\nabla\psi(x, s)|^2 e^{{i\psi(x,   s)}/\e}ds
+i\int^t_0  \hat{a}(x, s)c(s)\psi_s(x, s) e^{{i\psi(x,   s)}/\e}ds,&\\ \ns\ds
r_4=c(t)\big[\nabla\psi(x, t)\cdot \nabla\psi(x, t)-\psi_t^2(x, t)\big]e^{{i\psi(x, t)}/\e}.&
\end{array}\right.
\end{eqnarray}

 By \cite{FYZ,  M},  one  first may choose an $M(\cdot)\in
C^2([-T,  T];\;\dbC^{n\times n})$  with
$M(0)=M^0$,  such that
\begin{equation}\label{822-e1}
\nabla\psi(x, t)\cdot \nabla\psi(x, t)
-\psi_t^2(x, t)=\mathcal{O}(|x-\hat x(t)|^3),\hb{
as }x\to \hat x(t),\qq \forall\   t\in [-T, T].
\end{equation}
This implies that
\begin{equation}\label{xz1}
|r_4(\cd, t)|_{L^2(\dbR^n)}=\mathcal O(\e
^{\frac{n}{4}+\frac{3}{2}}),\;\; \hb{uniformly
for $\ae$} t\in (-T, T).
\end{equation}
Also,   one can choose  a $c(\cdot)\!\in C([-T,
T];\;\dbC\setminus\{0\})\bigcap
W^{2,\infty}((-T, T)\setminus\{0\})$ with
$c(0)=c^0$,  such that
$
2c_t(t)\psi_t(x,  t) +c(t)(W\psi)(x, t) =\mathcal
O(|x-\hat x(t)|),\hb{ as }x\to \hat x(t).
$
Meanwhile, one chooses  $A(s)=2i  c(s)  \hat{a}(\hat{x}(s), s)$.  Then
$$A(t)|\nabla\psi(x, t)|^2+i  c(t)\hat{a}(x,  t) \psi_t(x, t)=
|\hat a|_{L^\infty(0,  T;  W^{1,
\infty}(\dbR^n))}\mathcal{O}(|x-\hat{x}(t)|)\hb{
as }x\to \hat x(t),$$ and
\begin{equation}\label{xz8}
|r_3(\cd,t)|_{L^2(\dbR^n)}=(|\hat
a|_{L^\infty(0,  T;  W^{1, \infty}(\dbR^n))}+1)
\mathcal O(\e ^{\frac{n}{4}+\frac{1}{2}}),\q
\hb{uniformly for $\ae$} t\in (-T, T).
\end{equation}
Further, it is easy to check that
\begin{equation}\label{xz9}
|r_1(\cd,t)|_{L^2(\dbR^n)}=
|r_2(\cd,t)|_{L^2(\dbR^n)}=\mathcal O(\e
^{\frac{n}{4}}),\q \hb{uniformly for $\ae$} t\in
(-T, T).
\end{equation}
By \cite[Lemma 3.4]{FYZ}  and the definition of
$\phi_\e$,  for any $t\in [-T, T]$ and a
positive constant $c_2$,
\begin{equation}\label{**llz!}
|\phi_{\e}(\cdot,
t)|_{L^2(\dbR^n)}=\mathcal{O}(\e)  \mbox{ and }
|\phi_{\e, t}(\cdot, 0)|_{L^2(\dbR^n)},
|\phi_\e(\cdot, 0)|_{H^1_0(\dbR^n)}\leq c_2.
\end{equation}
Similar to arguments in \cite{FYZ}, by \eqref{Wu}-\eqref{**llz!}, one can easily get
the following results:

\vspace{0.1cm}

\noindent 1)  $ \{\phi_\e\}_{\e>0}$ given in $(\ref{ssrz5})$
is a sequence of  approximate solutions to \eqref{8.18-eq1} in
the sense that
$$\ds
\esssup_{t\in(-T, T)}\left| (W \phi_\e)(\cd,
t)+F(\phi_\e)(\cd,  t)
\right|_{L^{2}(\dbR^n)}=(|\hat a|_{L^\infty(0,
T;  W^{1, \infty}(\dbR^n))}+1)\mathcal
O(\e^{\frac{1}{2}}),\qq\hb{as }\e\to 0.
$$

\noindent 2) The initial energy of $\phi_\e$ is
bounded below, i.e.,
$\ds | \phi_{\e, t}(\cd, 0)|_{L^{2}(\dbR^n)}\ge
c_1 $
for a positive constant $c_1$, independent of
$\e$,  and
$\ds
|\phi_\e(\cd,  0)|_{L^2(\dbR^n)}=\mathcal
O(\e),\  \mbox{as }\e\to 0.
$

\medskip

\noindent 3) The energy of $\phi_\e$ is
polynomially  small off the generalized ray $(\hat x(\cdot),  p(\cdot))$:
$$\ds \esssup_{t\in (-T,
T)}\int_{\dbR^n\setminus B_{\e^{1/4}}(\hat
x(t))}\left(|\phi_{\e, t}(x,  t)|^2 +|
\phi_\e(x, t)|^2+|\n \phi_\e(x,
t)|^2\right)dx=\mathcal O(\e^2),\  \hb{as }\e\to 0.
$$
Here and hereafter, for any  $\k\subseteq
\dbR^n$
and $\d > 0$,   set
$
B_\d(\k)=\big\{ x\in \dbR^n\ \big|\  |x- x'| < \d \mbox{
for some } x' \in \k\big\}.
$

\smallskip

Furthermore,  we claim that
\begin{equation}\label{mmll}\ds
|\phi_{\e, t}(\cd, 0)|_{H^{-1}(\dbR^n)}=\mathcal
O(\e^{\frac{1}{2}}),\  \hb{as }\e\to 0.
\end{equation}
Indeed,  by (\ref{**llz!}),
$$
\ds|\nabla \phi_\e(\cdot,
0)|_{H^{-1}(\dbR^n)}\leq |\phi_\e(\cdot,
0)|_{L^2(\dbR^n)}= \mathcal{O}(\e).
$$
By \eqref{822-e1}, we get that
$$
\psi_t^2(x, 0) =\nabla\psi(x, 0)\cd
\nabla\psi(x, 0)+\mathcal{O}(|x-\hat
x(0)|^3),\qq\hb{ as }x\to \hat x(0).
$$
Hence,
\begin{equation}\label{822-e1.1}
|\psi_t(x,  0)|^2-
|\nabla\psi(x, 0)|^2=\mathcal{O}(|x-\hat
x(0)|^{3}),\qq\hb{ as }x\to \hat x(0).
\end{equation}
Therefore, by the definition of $\phi_\e$ and \cite[Lemma 3.4]{FYZ},
$$
\begin{array}{ll}
\ds |\phi_{\e,  t}(\cdot, 0)|_{H^{-1}(\dbR^n)}=
\big|\e^{1-\frac{n}{4}} c_t(0) e^{i\psi(\cdot,
0)/\e}+ \e^{-\frac{n}{4}} c(0)i\psi_t(\cdot, 0)
e^{i\psi(\cdot,  0)/\e}+ \e^{2-\frac{n}{4}} A(0)
e^{i\psi(\cdot,
0)/\e}\big|_{H^{-1}(\dbR^n)} \\
\ns\ds \leq\big|\e^{1-\frac{n}{4}} c_t(0)
e^{i\psi(\cdot,  0)/\e}\big|_{L^2(\dbR^n)}+
\big|\e^{-\frac{n}{4}} c(0)i\psi_t(\cdot, 0)
e^{i\psi(\cdot,  0)/\e}\big|_{H^{-1}(\dbR^n)}+
\big|\e^{2-\frac{n}{4}} A(0) e^{i\psi(\cdot,
0)/\e}\big|_{L^2(\dbR^n)} \\
\ns\ds \leq
\mathcal{O}(\e)+\Big|\e^{-\frac{n}{4}} \big|c(0)
\psi_t(\cdot, 0) e^{i\psi(\cdot,  0)/\e}\big|-
\big|\e^{-\frac{n}{4}} c(0)e^{{i\psi(\cd,
0)}/\e}\n\psi(\cd,0)\big|\Big|_{L^2(\dbR^n)}+
|\nabla
\phi_\e(\cdot, 0)|_{H^{-1}(\dbR^n)} \\
\ns\ds \leq
\mathcal{O}(\e)+\Big|\e^{-\frac{n}{4}} |c(0)|
\big(|\psi_t(\cdot, 0)|-\big|\nabla \psi(\cdot,
0)\big|\big) |e^{i\psi(\cdot,  0)/\e}|
\Big|_{L^2(\dbR^n)}.
\end{array}
$$
Now,   it is sufficient to estimate the last term in the above inequality.
Set
$$
\cA_\e=\Big\{x\in\dbR^n\ \Big|\  |\psi_t(x, 0)| +
\big|\nabla \psi(x,
0)\big|\leq\e\Big\}.
$$
Then
\begin{equation}\label{8.18-eq4}
\begin{array}{ll}\ds
\Big|\e^{-\frac{n}{4}} |c(0)| \big(|\psi_t(\cdot,
0)|-\big|\nabla \psi(\cdot, 0)\big|\big)
|e^{i\psi(\cdot,  0)/\e}| \Big|_{L^2(\dbR^n)}^2\\\ns\ds
= \e^{-\frac{n}{2}} |c(0)|^2
\int_{\cA_\e}\big(|\psi_t(x, 0)|-\big|\nabla
\psi(x, 0)\big|\big)^2 |e^{2i\psi(x,
0)/\e}|dx\\
\ns\ds\q + \e^{-\frac{n}{2}} |c(0)|^2
\int_{\dbR^n\setminus\cA_\e}\big(|\psi_t(x,
0)|-\big|\nabla \psi(x,
0)\big|\big)^2 |e^{2i\psi(x, 0)/\e}|dx\\
\ns\ds \leq \e^{2-\frac{n}{2}} |c(0)|^2
\int_{\cA_\e} |e^{2i\psi(x, 0)/\e}|dx\\\ns\ds
\quad +
\e^{-2-\frac{n}{2}} |c(0)|^2
\int_{\dbR^n\setminus\cA_\e}\big(|\psi_t(x,
0)|^2-\big|\nabla \psi(x,
0)\big|^2\big)^2 |e^{2i\psi(x, 0)/\e}|dx.
\end{array}
\end{equation}
It  is easy  to show that
\begin{equation}\label{8.18-eq5}
\e^{2-\frac{n}{2}} |c(0)|^2 \int_{\cA_\e}
|e^{2i\psi(x, 0)/\e}|dx \leq\e^{2-\frac{n}{2}}
|c(0)|^2 \int_{\dbR^n} |e^{2i\psi(x, 0)/\e}|dx
=\mathcal{O}(\e^2).
\end{equation}
Also, by \eqref{822-e1.1}, we find that
$$
(|\nabla\psi(x, 0)|^2
-|\psi_t(x,  0)|^2)^2=\mathcal{O}(|x-\hat
x(0)|^6),\q\hb{ as }x\to \hat x(0).
$$
Hence, by \cite[Lemma 3.4]{FYZ}  again,
\begin{equation}\label{8.18-eq6}
\begin{array}{ll}\ds
\e^{-2-\frac{n}{2}} |c(0)|^2
\int_{\dbR^n\setminus\cA_\e}\big(|\psi_t(x,
0)|^2-\big|\nabla \psi(x,
0)\big|^2\big)^2 |e^{2i\psi(x,
0)/\e}|dx\\
\ns\ds\leq  \e^{-2-\frac{n}{2}} |c(0)|^2
\int_{\dbR^n}\big(|\psi_t(x,
0)|^2-\big|\nabla \psi(x,
0)\big|^2\big)^2 |e^{2i\psi(x, 0)/\e}|dx
=\e^{-2-\frac{n}{2}}  \mathcal{O}(\e^{3+\frac{n}{2}}
)=\mathcal{O}(\e).
\end{array}
\end{equation}
By \eqref{8.18-eq4}-\eqref{8.18-eq6}, we get
the  desired  estimate \eqref{mmll}.

\smallskip

\noindent {\bf Step 2.}  We construct highly
concentrated approximate solutions  to    the
hyperbolic equation in a bounded domain.
Consider the following hyperbolic equation:
\begin{equation}\label{rrz4}
\left\{
\begin{array}{ll}\ds
\f_{tt}-\Delta \f+ F(\f)=0&\hb{ in }Q,\\
\ns\ds \f=0&\hb{ on }\Si.
\end{array}
\right.
\end{equation}

Assume that $(\hat x^-(\cdot), p^-(\cdot))$ is a
generalized  ray of  $W$ starting from $\hat x
^-(0)\in \Omega$ and arriving  $ x_0= \hat
x^-(t_0)\in
\partial \Omega$. By Step 1, one can construct
a family of approximate solutions
$\{\phi_\e^{-}\}_{\e>0}$ to    the first equation
of (\ref{rrz4}). However, $\phi^-_{\e}$ may not
satisfy the homogeneous Dirichlet boundary
condition on $\Sigma$. To solve this problem, we
superpose $ \phi^-_{\e}$ with another
approximate solution $ \phi^+_{\e}$. The latter
is constructed from a  ray $(\hat x^+(\cdot),
p^+(\cdot))$, which reflects $(\hat x^-(\cdot),
p^-(\cdot))$ at the boundary $\partial \Omega$.
The  key point is to select an approximate
solution $\phi^+_{\e}$  concentrated in a small
neighborhood of the reflected ray $(\hat
x^+(\cdot), p^+(\cdot))$, such that
$\phi^-_{\e}+\phi^+_{\e}$ satisfies
approximately the homogeneous Dirichlet boundary
condition.

%According to (\ref{rz3}), $( \hat x^-(t),  \xi^-(t))$
%satisfies
%\begin{equation}\label{1z3}
%\begin{cases}
%\ds x_t^-(t)=2\xi^-(t),\\
%\ns \ds
%\xi_t^-(t)=0,\\
%\ns \ds \hat x^-( t_0)= x_0,\q
%\xi^-( t_0)= \xi^-( t_0).
%\end{cases}
%\end{equation}
%

\smallskip

Choose $( \hat x^+(\cdot),  p^+(\cdot))$,  such that
\begin{equation*}\label{0001z4}\left\{
\begin{array}{ll}
\ds \hat x_t^+(t)=-2p^+(t),\   p_t^+(t)=0,&\\ \ns
\ds \hat x^+( t_0)= x_0,\  p^+( t_0)=p^-( t_0)-2
\nu^\top(x_0) p^-( t_0) \nu( x_0).
\end{array}\right.
\end{equation*}
%
%
%It follows from (\ref{OK11}) that
%$|\xi^-(t_0)|^2=\frac{1}{4}$. Hence,
%$|\xi^+(t_0)|^2=\frac{1}{4}$. From the second
%equation of \eqref{1z3} and \eqref{0001z4}, we
%see that
%
%Then it is easy to show that
%\begin{equation}\label{KO22}
%|p^\pm(t)|^2
%=\frac{1}{4},\qq\forall\;t\in\dbR.
%\end{equation}
%
Assume that $p^-(\cdot)$ is transversal to the
boundary $\partial \Omega$ at  the time $t_0$, i.e., $
[p^-(t_0)]^\top \nu(x_0)\neq 0. $ Denote by
$t_1>0$ the instant,  when the reflected ray
arrives at $\partial \Omega$, i.e., $ \hat x^+(t_1)\in \partial \Omega$.
For any $ T^*\in ( t_0, t_1)$, choose a cut-off
function $\varrho^-\in C_0^\infty(\dbR^{n+1})$,
which  equals to   $1$ identically in a
neighborhood of the set $\big\{(t,\hat
x^-(t))\;\big|\;t\in [0,t_0]\big\}$ with $
 \supp\varrho^-\subseteq  B_{(T^*-t_0)/4}\big\{(t,\hat x^-(t))\;\big|\;t\in
[0, t_0]\big\}.
$   Then
by Step 1, we may construct approximate
solutions to \eqref{8.18-eq1} as
follows:
%
%\begin{equation}\label{000rz5}
$$
\phi_\e^{-}(x,  t)=\e^{1-n/4}\varrho^-(x,   t)
c^-(t)e^{i \psi^-(x, t)/\e },
$$
%\end{equation}
%
where
%\begin{equation}\label{1-rz6}
$$
\psi^-(x, t)=[p^-(t)]^\top[x- \hat
x^-(t)]+\frac{1}{2}\big[x- \hat x^-(t)\big]^\top
M^-(t)\big[x- \hat x^-(t)\big]
$$
%end{equation}
and  $M^-(t)$ is a
complex symmetric matrix with positive definite
imaginary part.

\smallskip

Next, we  construct another approximate solution
to  \eqref{8.18-eq1}  as follows:
\begin{equation}\label{0000rz5}
\phi_\e^{+}(x,  t)=\e^{1-\frac{n}{4}}\varrho^+(x,  t)c^+(t)e^{i
\psi^+(x,  t)/\e },
\end{equation}
 which is
concentrated in a  small  neighborhood of the
reflected ray  $(\hat x^+(\cdot), p^+(\cdot))$, such
that
\begin{equation}\label{1z7}
|\phi_\e^-+\phi_\e^{+}|_{H^1(\partial \Omega\times(0,T^*))}=\mathcal O(\e^{1/2}).
\end{equation}
In \eqref{0000rz5}, $\varrho^+\in
C_0^\infty(\dbR^{n+1})$ is a  cut-off function,
which identically equals to $1$ in a
neighborhood of the set $\big\{(t,\hat
x^+(t))\;\big|\;t\in [t_0,T_1]\big\}$,  such
that
%
%\begin{equation}\label{or21}
$
\supp\varrho^+\subseteq
B_{\min\{t_0,T_1-T^*\}/4} \big\{(t,\hat
x^+(t))\;\big|\;t\in [t_0,T_1]\big\},
$
%\end{equation}
and $$ \psi^+(t,x)= [p^+(t)]^\top[x- \hat
x^+(t)]+\frac{1}{2}\big[x- \hat x^+(t)\big]^\top
M^+(t)\big[x- \hat x^+(t)\big].
$$
Here $M^+(\cdot)$ and $c^+(\cdot)$ are  determined in
a similar way in Step 1. The only difference is
that $
c^+(t_0)=- c^-(t_0).
$
Also,  $M^+(t)$ is determined by its initial
$M^+(t_0)$ and the reflected ray $(\hat
x^+(\cdot), p^+(\cdot))$. Similar to the choices in
\cite{FYZ,   M}, it is easy to check
that \eqref{1z7} holds.

\smallskip

\noindent {\bf Step 3.}  We construct
approximate solutions $\{\Phi_\e\}_{\e>0}$ to
\eqref{rrz4},  such that the energies of them
are concentrated in a neighborhood of  a generalized ray
$\big\{(t, \hat{x}^j(t), p^j(t))\;\big|\;
t\in[s_j,s_{j+1}]\big\}_{j=0}^{m-1}$. To this end,
choose a cut-off function $\varrho^1\in
C_0^\infty(\dbR^{n+1})$,  which identically
equals  to $1$ in a neighborhood of
$\big\{(t,\hat x^1(t))\;\big|\;t\in
[0,s_1]\big\}$ with
$\supp\varrho^1\subseteq B_{(s_2-s_1)/4}
\big\{(t,\hat x^1(t))\;\big|\;t\in
[0,s_1]\big\}. $
By Step 2, we can find a function
$ \phi_\e^1(x,  t)=\e^{1-n/4}\varrho^1(x, t)
c_1(t)e^{i \psi^1(x,  t)/\e },$
such that
\begin{equation*}\label{or7}\left\{
\begin{array}{ll}\ds
\esssup_{t\in(0,s_1)}\left|(W\phi_\e^1)(\cd, t)+ F(
\phi_\e^1)(\cd, t)\right|_{L^{2}(\O)}=(|\hat a|_{L^\infty(0,  T;  W^{1, \infty}(\Omega))}+1)\mathcal O(\e^{\frac{1}{2}}),&\\
\ns\ds |\phi_{\e, t}^1(\cd, 0)|_{L^2(\O)}\ge
c_1,\  |\phi^1_{\e, t}(\cdot, 0)|_{L^2(\O)},   |\phi^1_\e(\cdot, 0)|_{H^1_0(\O)}\leq  c_2,&\\
\ns\ds
|\phi^1_\e(\cdot, 0)|_{L^2(\Omega)}+|\phi^1_{\e, t}(\cdot, 0)|_{H^{-1}(\Omega)}=
\mathcal{O}(\e^{1/2}),&\\ \ns\ds
|\phi_\e^1|_{H^1(0,T;\,L^2(\o))}+|\phi_\e^1|_{H^1(-T,0;\,L^2(\o))}=\mathcal O(\e^{\frac{1}{2}}).
\end{array}\right.
\end{equation*}

Next, choose a cut-off function $\varrho^2\in
C_0^\infty(\dbR^{n+1})$,  which identically
equals to $1$ in a neighborhood of
$\big\{(t,\hat x^2(t))\;\big|\;t\in
[s_1,s_2]\big\}$ with $ \supp\varrho^2\subseteq
B_{\min\{s_1,s_3-s_2\}/4}\big\{(t,\hat
x^2(t))\;\big|\;t\in [s_1,s_2]\big\}. $ By  Step
2, we can find a  $ \phi_\e^2(x,
t)=\e^{1-n/4}\varrho^2(x, t) c_2(t)e^{i
\psi^2(x,  t)/\e }, $ such that for  $S\in
(0,s_2)$,
\begin{equation*}\label{or9}\left\{
\begin{array}{ll}\ds
\esssup_{t\in(s_1,s_2)}\left|(W \phi_\e^2)(\cd,  t)+F(\phi_\e^2)(\cd,  t)\right|_{L^2(\O)}=
(|\hat a|_{L^\infty(0,  T;  W^{1, \infty}(\Omega))}+1)\mathcal O(\e^{\frac{1}{2}}),&\\
\ns \ds |\phi_\e^1+\phi_\e^2|_{H^1(\partial \Omega\times(0,S))}+
|\phi_\e^1+\phi_\e^2|_{H^1(\partial \Omega\times(-T,0))}=\mathcal O(\e^{\frac{1}{2}}),&\\
\ns\ds
|\phi_\e^2|_{H^1(0,T;\,L^2(\o))}+|\phi_\e^2|_{H^1(-T,0;\,L^2(\o))}=\mathcal O(\e^{\frac{1}{2}}).
\end{array}\right.
\end{equation*}
Further, for $j=3,\cdots, m$, choose a cut-off function
$\varrho^j\in C_0^\infty(\dbR^{n+1})$, which
identically equals to $1$ in a neighborhood of
$\big\{(t,\hat x^j(t))\;\big|\;t\in
[s_{j-1},s_j]\big\}$ with
$
\supp\varrho^j\subseteq B_{\min\{s_{j-1}-s_{j-2},s_{j+1}-s_j\}/4}
\big\{(t,x^j(t))\;\big|\;t\in
[s_{j-1},s_j]\big\}.
$
Similarly, we can find a
$
\phi_\e^j(x, t)=\e^{1-n/4}\varrho^j(x,  t)
c^j(t)e^{i\psi^j(x,  t)/\e },
$
such that
\begin{equation*}\label{or10}\left\{
\begin{array}{ll}\ds
\esssup_{t\in(s_{j-1},s_j)}\left|(W
\phi_\e^j)(\cd,  t)+F(\phi_\e^j)(\cd,  t)\right|_{L^2(\O)}=(|\hat a|_{L^\infty(0,  T;  W^{1, \infty}(\Omega))}+1)\mathcal O(\e^{\frac{1}{2}}),&\\
\ns
\ds |\phi_\e^{j-1}+\phi_\e^j|_{H^1(\partial \Omega\times(0, S))}+|\phi_\e^{j-1}+\phi_\e^j|_{
H^1(\partial \Omega\times(-T, 0))}=\mathcal O(\e^{\frac{1}{2}}),&\\
\ns\ds
|\phi_\e^j|_{H^1(0,T;\,L^2(\o))}+|\phi_\e^j|_{H^1(-T,0;\,L^2(\o))}=\mathcal O(\e^{\frac{1}{2}}),
\end{array}\right.
\end{equation*}
for $S\in (0,s_j),  \mbox{ if  }
j=3,\cdots,m-1,\mbox{ and }S\in (0,T], \mbox{ if
}j=m.$

Now,  write
$
\Phi_\e=\sum\limits_{j=1}^m \phi_\e^j.
$
Then it is easy to show that \begin{equation*}\left\{
\begin{array}{ll}\ds
\esssup_{t\in(0,T)}\left|(W
\Phi_\e)(\cd,  t)+F(\Phi_\e)(\cd,   t) \right|_{L^2(\O)}=(|\hat a|_{L^\infty(0,  T;  W^{1, \infty}(\Omega))}+1)\mathcal O(\e^{\frac{1}{2}}),&\\
\ns
\ds  |\Phi_\e|_{H^1(\partial \Omega\times(0,  T))} + |\Phi_\e|_{H^1(\partial \Omega\times(-T, 0))}=\mathcal O(\e^{\frac{1}{2}}),&\\
\ns\ds
|\Phi_{\e,  t}(\cd,  0)|_{L^2(\O)}\ge c_1, \  |\Phi_{\e, t}(\cdot, 0)|_{L^2(\O)},   |\Phi_\e(\cdot, 0)|_{H^1_0(\O)}\leq  c_2,&\\
\ns\ds |\Phi_\e(\cd, 0)|_{L^2(\O)} + |
\Phi_{\e,  t}(\cd,  0)|_{H^{-1}(\O)}=\mathcal O(\e^{\frac{1}{2}}),&
\\
\ns\ds
|\Phi_\e|_{H^1(0,T;\,L^2(\o))}+|\Phi_\e|_{H^{1}(-T,0;\,L^2(\o))}=\mathcal  O(\e^{\frac{1}{2}}).
\end{array}\right.
\end{equation*}
Set
$
\wt\Phi_\e(x, t)=\Phi_\e(x,  t)+\Phi_\e(x, -t),
\mbox{ for } (x, t)\in Q.
$
Then
$\{\wt\Phi_\e\}_{\e>0}$ satisfies
\begin{equation}\label{or13.2}\left\{
\begin{array}{ll}\ds
\esssup_{t\in(0,T)}\left|(W
\wt\Phi_\e)(\cd, t)+F(\wt\Phi_\e)(\cd,  t) \right|_{L^2(\O)}=(|\hat a|_{L^\infty(0,  T;  W^{1, \infty}(\Omega))}+1)\mathcal O(\e^{\frac{1}{2}}),&\\
\ns
\ds  |\wt\Phi_\e|_{H^1(\partial \Omega\times (0, T))} =\mathcal O(\e^{\frac{1}{2}}),&\\
\ns\ds |\wt\Phi_{\e, t}(\cd, 0)|_{L^2(\O)}\ge
c_1,\  |\wt\Phi_{\e, t}(\cdot, 0)|_{L^2(\O)}\ |\wt\Phi_\e(\cdot, 0)|_{H^1_0(\O)}\leq  c_2,& \\
\ns\ds |\wt\Phi_\e(\cd,  0)|_{L^2(\O)} +|
\wt\Phi_{\e, t}(\cd, 0)|_{H^{-1}(\O)}=\mathcal  O(\e^{\frac{1}{2}}),&
\\
\ns\ds |\wt\Phi_\e|_{H^1(0,T;\,L^2(\o))}
=\mathcal O(\e^{\frac{1}{2}}).
\end{array}\right.
\end{equation}

\smallskip

\noindent {\bf Step 4.} We construct a family of
solutions $\{\phi_\e\}_{\e>0}$ to  \eqref{rrz4}.  To this aim, let
$
\phi_\e=\wt\Phi_\e+v_\e,
$
where $v_{\e}$ solves
\begin{equation*}\label{or14}
\left\{
\begin{array}{ll}\ds
Wv_\e  + F(v_\e) =-W\wt\Phi_\e - F(\wt\Phi_\e) &\hb{ in }Q,\\
\ns\ds
v_\e=-\wt\Phi_\e\qq\qq\qq\qq &\hb{ on }\Si,\\
\ns\ds v_{\e}(x, 0)=0,\   v_{\e, t}(x, 0)=0 \qq&\hb{
in } \O.
\end{array}
\right.
\end{equation*}
It is easy to see that
$$
\begin{array}{ll}\ds
\max_{t\in[0,T]}|(v_\e(\cd, t),
v_{\e, t}(\cd,  t))|_{H^{1}(\O)\times L^2(\O)}\ds\le C\big(|W\wt\Phi_\e +
F(\wt\Phi_\e)|_{L^1(0,T;\,L^2(\O))}+|\wt\Phi_\e|_{H^1(\partial \Omega\times(0, T))}\big).
\end{array}
$$
This, together with the first two conclusions in
\eqref{or13.2},  implies that
\begin{equation}\label{or16}
\max_{t\in[0,T]}|(v_\e(\cdot, t),
v_{\e,  t}(\cd,  t))|_{H^{1}(\O)\times L^2(\O)}\le
C(|\hat a|_{L^\infty(0,  T;  W^{1, \infty}(\Omega))}+1)\e^{\frac{1}{2}}.
\end{equation}
By   \eqref{or13.2} and \eqref{or16}, we get
that
\begin{equation*}\label{5.8-eq1}\left\{
\begin{array}{ll}\ds
|\phi_{\e, t}(\cd,  0)|_{L^2(\O)}\ge c_1, \  |\phi_{\e, t}(\cdot, 0)|_{L^2(\O)}, |\phi_\e(\cdot, 0)|_{H^1_0(\O)}\leq c_2,&\\
\ns\ds |\phi_\e(\cd, 0)|_{L^2(\O)} + |
\phi_{\e,  t}(\cd,  0)|_{H^{-1}(\O)}=\mathcal O(\e^{\frac{1}{2}}),&
\\ \ns\ds |\phi_\e|_{H^1(0,T;\,L^2(\o))}
=(|\hat a|_{L^\infty(0,  T;  W^{1, \infty}(\Omega))}+1)\mathcal  O(\e^{\frac{1}{2}}).
\end{array}\right.
\end{equation*}
Hence,
$\{\phi_\e\}_{\e>0}$ are the desired family of solutions in   Lemma \ref{8.18-lm1}.
\endpf

\end{document}